\documentclass [a4paper,12pt]{article}
\usepackage[dvips]{graphicx}
\usepackage[english]{babel}
\usepackage{amssymb,latexsym}
\usepackage{mathrsfs}
\usepackage[all]{xy}
\usepackage{color}

\usepackage{epsfig,afterpage}
\usepackage{psfrag}
\usepackage{graphicx}

\usepackage[ansinew]{inputenc}
\input epsf

\usepackage{amsmath}
\usepackage{amssymb}
\usepackage{amsfonts}
\usepackage{multirow}

\newtheorem{theorem}{Theorem}
\newtheorem{lemma}{Lemma}
\newtheorem{remark}{Remark}

\newtheorem{exam}{Example}
\newtheorem{defin}{Definition}

\newcommand{\ebox}{\fbox {} \smallskip}

\numberwithin{equation}{section}  

\def\{{\protect\lbrace}
\def\}{\protect\rbrace}
\def\tr{\operatorname{tr}}
\def\Real{\operatorname{Re}}

\def\const{{\rm{const}}}
\newcommand{\I}{\bigl |}

\newcommand{\IIII}{\Biggl |}
\newcommand{\pa}{\partial}

\def\alf{\alpha}
\def\bet{\beta}
\def\Del{\Delta}
\def\del{\delta}
\def\eps{\varepsilon}
\def\gam{\gamma}
\def\Gam{\Gamma}
\def\phi{\varphi}

\def\om{\omega}
\def\la{\lambda}

\def\ov{\overline}

\newcommand{\bR}{\mathbb{R}}

\newcommand{\bN}{\mathbb{N}}
\newcommand{\bI}{\mathbb{I}}

\newcommand{\FF}{\mathscr {F}}
\newcommand{\F}{\Phi}
\newcommand{\MM}{{\mathfrak M}}

\def\oF{{\ov F}}
\def\oH{\ov H}
\def\oG{\ov G}
\def\of{{\ov f}}

\renewcommand{\AA}{A}
\newcommand{\FFF}{\mathcal B}

\def\dx{{\dot x}}
\def\dy{{\dot y}}
\def\ddx{{\ddot x}}
\def\ddy{{\ddot y}}
\def\id{{\rm id}}
\def\M{M}

\let\<\langle
\let\>\rangle

\title
{
Geodesics in generalized Finsler spaces: \\
singularities in dimension two
}

\begin{document}

\author{
A.O. Remizov\footnote{Supported by the FAPESP visiting professor grant 2014/04944-6.}
}

\maketitle

\begin{abstract}
We study singularities of geodesics flows in two-dimensional generalized Finsler spaces (pseudo-Finsler spaces).
Geodesics are defined as extremals of a certain auxiliary functional whose non-isotropic extremals coincide with
extremals of the action functional. This allows to consider isotropic lines as (unparametrized) geodesics.
\end{abstract}

\renewcommand{\thefootnote}{\fnsymbol{footnote}}

\footnote[0]{Instituto de Ci\^encias Matem\'aticas e de Computa\c{c}\~ao - USP,
Avenida Trabalhador S\~ao-Carlense, 400 - Centro,
CEP: 13566-590 - S\~ao Carlos - SP, Brazil.
Email: alexey-remizov@yandex.ru
}
\footnote[0]{2010 Mathematics Subject classification 53B40 (Primary), 53C22, 34C05.}

\section*{Introduction}

The paper is a study of singularities of geodesics flows in generalized Finsler spaces
(pseudo-Finsler spaces). This is a natural development of an ongoing research on understanding the geometry
of surfaces endowed with a signature varying pseudo-Riemannian metric;
see \cite{GR, KhT, KosKri, PR, Rem-Pseudo, Rem15, RemizovTari} and the references therein.
One of the purposes of the paper is to compare singularities of geodesics flows
in pseudo-Finsler and pseudo-Riemannian metrics.
On the other hand, the interest to pseudo-Finsler spaces is motivated by physical applications;
see e.g. \cite{AM, Asanov}.

According to \cite{Rund}, by a pseudo-Finsler space we mean a manifold $M$, $\dim M = m$, with coordinates $(x_i)$
endowed with a metric function $\of(x_i;\dx_i) = \oF(x_i;\dx_i)^{\frac{1}{n}}$,
where $\oF \colon TM \to \bR$ is positively homogeneous in $(\dx_i)$ of degree
$n$ and smooth on the complement of the zero section of $TM$
(a more detailed definition is given in Sections \ref{sect1.1}, \ref{sect1.2}).
A well-known example is {\it Berwald--Moor} space $(M,\of)$, where
$\of(x_i;\dx_i) = (\dx_1 \cdots \dx_n)^{\frac{1}{n}}$, $n=m$; see, e.g., \cite{Balan, Matsumoto-Shimada, Rund}.

The paper starts with the discussion of the notion of geodesics in Finsler and pseudo-Finsler spaces
with $n \ge 3$ (Section~\ref{sect1}). Here we use the variational definitions of geodesics \cite{BCS, Rund}.
In contrast to pseudo-Riemannian spaces ($n=2$), where naturally pa\-ra\-met\-rized geodesics of all types (including isotropic) can be defined as extremals of the action functional,
in pseudo-Finsler spaces the similar definition is not correct for isotropic lines.
The solution of this problem is either to exclude isotropic lines from consideration or to find a natural extension
of the definition of geodesics.

In the present paper, we choose the second way.
Based on a simple variational property, we define geodesics as extremals of a certain auxiliary functional whose non-isotropic extremals coincide with extremals of the action functional.
In this direction, we have the following result: in the case $m=2$, all isotropic lines are (unparametrized) geodesics.

In Section \ref{sect2}, we consider singularities of the geodesic flows in pseudo-Finsler spaces $(M, \of)$,
where $m=2$ and $\oF$ is a polynomial in $(\dx_i)$ of degree $n \ge 3$.
The main results are presented in Section~\ref{sect2.2}, where we consider the case $n=3$ in detail.
It is proved that singularities of the geodesic flow are connected with the degeneracy of isotropic lines net.
Namely, if the function $\oF$ is generic, the manifold $M$ contains two open domains $M_+$ and $M_-$
separated by a curve $M_0$ so that at every point $q \in M_+$ (resp., $M_-$) there exist 3 (resp., 1) different
isotropic directions and isotropic lines are tangent at $q \in M_0$.
Singularities of the geodesic flow appear in the domain $M_-$ and on the curve $M_0$.

Section \ref{sect3} is devoted to a special case: pseudo-Finsler space $(M, \of)$, where $M$ is a surface
in $n$-dimensional Berwald--Moor space. The corresponding function $\oF$ is a non-generic
polynomial in $(\dx_i)$ of degree $n$. In this case, the domain $M_- = \emptyset$, and
singularities of the geodesic flow appear on the curve $M_0$ only.

\medskip

The author expresses deep gratitude to prof. Farid Tari (ICMC-USP, S\~ao Carlos, Brazil)
for attention to the work, useful advices, and remarks.
I am also grateful to the reviewer for many constructive comments and suggestions.

\section{Variational definition of geodesics}\label{sect1}

\subsection{Finsler spaces}\label{sect1.1}

Consider a smooth
(here and below, by smooth we mean $C^{\infty}$ unless otherwise stated)
manifold $M$, $\dim M = m$, with coordinates $(x_i)$ and
a function $\oF(x_i;\dx_i) \colon TM \to \bR$ that is
positively homogeneous of degree $n \ge 2$ in $(\dx_i)$ and
smooth on the complement of the zero section of $TM$.

Define the function $\of(x_i;\dx_i) = \oF^{\frac{1}{n}}(x_i;\dx_i)$, which is
positively homogeneous of degree $1$ in $(\dx_i)$.
The pair $(M, \of)$ or, equivalently, $(M, \oF)$ is a {\it Finsler space} (in the classic sense),
if the following conditions hold:
\begin{itemize}
\item[{\bf B}.]
$\oF(x_i;\dx_i) > 0$ if $|\dx_1|+ \cdots+ |\dx_m| \neq 0$.
\item[{\bf C}.]
The Hessian of the function $\of^2$ with respect to $(\dx_i)$ is positive definite, that is,
\begin{equation}
\sum_{i,j=1}^m \frac{\pa^2 \bigl({\of}^2\bigr)}{\pa \dx_i \pa \dx_j}\, \xi_i \xi_j > 0
\ \ \, \textrm{if} \ \ \, \sum_{i=1}^m |\xi_i| \neq 0.
\label{1}
\end{equation}
\end{itemize}
Here we use the letters B and C to preserve the notations from the book \cite{Rund},
to which we shall refer.
The quadratic form \eqref{1} is called {\it the fundamental tensor}, and $\of$
(positive and smooth on the complement of the zero section of $TM$) is called the {\it metric function} on $M$.

The metric function $\of (x_i;\cdot)$ defines a Minkowski norm on each tangent space $T_xM$.
For a curve $\gamma \colon I \to M$, it allows to
define the length and the action functionals similarly to the Riemannian metrics:
\begin{equation}
J^{(\nu)}(\gamma) = \int\limits_{I} {\of}^{\,\nu}(x_i;\dx_i) \,dt =
\int\limits_{I} {\oF}^{\,\frac{\nu}{n}}(x_i;\dx_i) \,dt, \ \ \,
\dx_i = \frac{dx_i}{dt},
\label{2}
\end{equation}
with $\nu = 1$ (length) and $\nu = 2$ (action), see, e.g., \cite{Matsumoto95, Rund}.
As in the Riemannian  case, the length functional $J^{(1)}$ is invariant with respect to reparametrizations of $\gamma$, while $J^{(2)}$ is not.

Parametrized geodesics can be defined as extremals of the action functional $J^{(2)}$,
the corresponding parametrization is called {\it natural} or {\it canonical}
(it coincides with the arc-length parametrization, where $ds=\of$).

Non-parametrized geodesics can be defined as extremals of any one of the functionals $J^{(2)}$ and $J^{(1)}$.
The difference between using $J^{(2)}$ and $J^{(1)}$ is the following.
In the first case we just simply forget the natural parametrization of the extremals of $J^{(2)}$, while
in the second case the Euler-Lagrange system with the Lagrangian $\of (x_i;\dx_i)$ contains $m-1$
independent equations only \cite{Rund}. This reflects the fact that the length functional $J^{(1)}(\gamma)$
is invariant with respect to reparametrizations of $\gamma$.
Using this degree of freedom and assuming that we deal with continuously differentiable geodesics
with definite tangent directions at all points, one can put (at least, locally)
the parameter $t$ equal to one of the coordinates $x_i$,
and consequently, reduce the Euler-Lagrange system for $J^{(1)}(\gamma)$.

From now on, we use the following general notation.
Let $\ov \Phi(x_i;\dx_i)$ be a function on $T\M$ positively homogeneous of degree $k$ in $(\dx_i)$,
then the formula
\begin{equation}
\Phi = \frac{\ov \Phi}{\dx_1^k}
\label{def}
\end{equation}
defines a function on the projectivized tangent bundle $PT\M$.
For instance, put $x=x_1$ and $y_i=x_i$, $p_i=dy_i/dx$ for $i=2, \ldots, m$.
This yields
\begin{equation}
\frac{d}{dx} \biggl( \frac{\pa f}{\pa p_i} \biggr) = \frac{\pa f}{\pa y_i}, \ \ \
f(x,y_i,p_i)= \frac{\of (x_i;\dx_i)}{\dx_1}, \ \ \ i=2, \ldots, m.
\label{3}
\end{equation}
The passage to equation \eqref{3} is the standard projectivization
$\Pi \colon T\M \to PT\M$ of the tangent bundle.
Moreover, non-parametrized geodesics can be defined as extremals of $J^{(\nu)}$
with arbitrary $\nu \ge 1$, on the basis of the following simple property (see, e.g., \cite{MHV}).

\begin{lemma}
\label{L1}
Let $\psi \colon \bR \to \bR$ be a function such that $\psi \circ {\oF}$ is $C^2$-smooth
on the complement of the zero section of  $TM$ and $\psi'(s) \neq 0$ for all $s \neq 0$.
Then non-parametrized extremals of the functional
\begin{equation}
J_{\psi}(\gamma) = \int\limits_{I} \psi \circ {\oF}(x_i;\dx_i) \,dt, \ \ \,
\dot x = \frac{dx}{dt}, \ \ \, \dot y = \frac{dy}{dt},
\label{4}
\end{equation}
coincide with non-parametrized extremals of $J_{\id}(\gamma)$, where $\id$ is the identity map.
\end{lemma}

{\bf Proof}.
The Euler-Lagrange equation of $J_{\psi}(\gamma)$ reads
\begin{equation}
\frac{d}{dt} \biggl( \psi' \circ {\oF}(x_i;\dx_i)\cdot \frac{\pa \oF}{\pa \dx_i} \biggr) =
\psi' \circ {\oF}(x_i;\dx_i)\cdot \frac{\pa \oF}{\pa x_i}, \ \ \ i=1, \ldots, m.
\label{5}
\end{equation}
In light of the condition $\oF(x_i;\dx_i) > 0$, every curve $\gam$ admits the arc-length
parametrization, that is, ${\oF}(x_i;\dx_i) \equiv c \neq 0$ along $\gam$.
Using the arc-length parametrization in \eqref{4},
after reducing the constant factor $\psi' \circ {\oF}(x_i;\dx_i) = \psi'(c)$ in both sides of \eqref{5},
we get the Euler-Lagrange equation of the functional $J_{\id}(\gamma)$.
\hfill\ebox

Thus non-parametrized geodesics can be defined as extremals of $J_{\psi}(\gamma)$
with arbitrary function $\psi$ from Lemma~\ref{L1}, in particular, of the functional $J_{\id}(\gamma)$,
which is equal to $J^{(\nu)}(\gamma)$ with $\nu=n$ from \eqref{2}.
In classical Finsler spaces this extended definition of geodesics gives us nothing essentially new,
but it can be useful for generalized Finsler spaces considered below.

\subsection{Generalized Finsler spaces}\label{sect1.2}

A generalization of the notion of a Finsler space may be obtained if conditions B and C are dropped,
such spaces are sometimes called {\it special Finsler} or {\it pseudo-Finsler}.
Here we take the liberty to cite a passage from the classical book \cite{Rund} (page 265):
\vspace{-1.0ex}
\begin{quote}
{\it
Again, it should be remarked that very frequently the metric function
which is given by a homogeneous Lagrangian function of a dynamical
system does not always satisfy conditions B and C.
The singularities which may occur as a result of the relaxation of condition C
are usually ignored, but it is well possible that an investigation of these
singularities in connection with physical applications cannot be avoided
and might furthermore prove to be fruitful.
}
\end{quote}
\vspace{-1.0ex}

From now on, we will consider pseudo-Finsler spaces $(M, \oF)$,
where the function $\oF$ is not supposed to satisfy the conditions B and C.

The absence of the condition B brings to the existence of the {\it isotropic hypersurface} $\FF$
given by the equation $\oF (x_i;\dx_i) = 0$ in $TM$ or, equivalently, $F(x,y_i,p_i)=0$ in $PTM$.
The Euler-Lagrange equation for the functional $J^{(\nu)}(\gamma)$ with $\nu < n$ is not defined on $\FF$,
since the derivatives of ${\oF}^{\,\frac{\nu}{n}}(x_i;\dx_i)$ are discontinuous on $\FF$.
This explains the advantage of the functional $J^{(n)}(\gamma)$ for definition of geodesics
compared to $J^{(\nu)}(\gamma)$ with $\nu < n$.

The Euler-Lagrange equation for the functional $J^{(n)}(\gamma)$ reads
\begin{equation}
\sum_{j=1}^{m} \frac{\pa^2 \oF}{\pa \dx_i \pa \dx_j} \, \ddx_j +
\sum_{j=1}^{m} \frac{\pa^2 \oF}{\pa \dx_i \pa x_j} \, \dx_j =
\frac{\pa \oF}{\pa x_i}, \ \ \ i=1, \ldots, m,
\label{6}
\end{equation}
or, equivalently,
\begin{equation}
\ddx_i = \frac{\oH_i(x_i;\dx_i)}{\oH(x_i;\dx_i)}, \ \ \, \textrm{where} \ \ \,
\oH = \det \biggl( \frac{\pa^2 \oF}{\pa \dx_i \pa \dx_j} \biggr), \ \ \ i=1, \ldots, m,
\label{7}
\end{equation}
and $\oH_i$ are the determinants defined from \eqref{6} by Cramer's rule.
It is not hard to see that the functions $\oH_i$ are positively homogeneous of degree $n$ in $(\dx_i)$
and $\oH$ is positively homogeneous of degree $n-2$ in $(\dx_i)$.

Similarly to \eqref{3}, the projectivization $\Pi \colon T\M \to PT\M$ sends equation \eqref{7} to
\begin{multline}
p_i = \frac{dy_i}{dx}, \ \ \
\frac{dp_i}{dx} = \frac{dp_i}{dt} \biggl(\frac{dx}{dt}\biggr)^{-1} \!\!\! =
\frac{1}{\dx} \frac{d}{dt} \biggl(\frac{\dy_i}{\dx}\biggr) =
\frac{\ddy_i \dx - \dy_i\ddx}{\dx^3} = \\
\frac{1}{\dx^2} \frac{\oH_i - p_i\oH_1}{\oH} = \frac{H_i - p_iH_1}{H}, \ \ \ i=2, \ldots, m,
\label{8}
\end{multline}
where the functions $H, H_i$ are obtained from $\ov H, \ov H_i$ by \eqref{def}.
Recall that, by formula \eqref{def}, the independent variable $x$ is the coordinate $x_1$.
From Lemma \ref{L1}, it follows that out of the isotropic hypersurface $\FF$,
integral curves of \eqref{8} coincide with integral curves of \eqref{3}.
However, equation \eqref{8} is defined in the whole space $PTM$.

\begin{lemma}
\label{L2}
$\FF$ is an invariant hypersurface of both equations \eqref{6} and \eqref{8}.
Moreover, in the case $m=2$ all isotropic lines are non-parametrized extremals
of the functional $J^{(n)}$.\footnote{
This statement holds true only for $m=2$, while in the case $m>2$ some isotropic lines are geodesics
and some of them are not. An example is given for $m=3$, $n=2$ (pseudo-Riemannian metrics) in \cite{Rem15}.
}
\end{lemma}

{\bf Proof}.
After straightforward transformations, equation \eqref{3}
gives a direction field which is parallel to
\begin{equation}
F^{-\mu} X \biggl( \frac{\pa}{\pa x} +  p \frac{\pa}{\pa y} \biggr) +
F^{-\mu} \sum_{i=2}^{m} Y_i \frac{\pa}{\pa p_i},
\label{9}
\end{equation}
where $\mu = (m-1)(2-\frac{1}{n})$ and $X, Y_i$ are smooth functions on $PTM$.

Since the vector field \eqref{9} is derived directly from the Euler-Lagrange equation \eqref{3},
it is divergence-free in $PTM$ except for the hypersurface $\FF$ where the factor $F^{-\mu}$ is discontinuous, and the field is not defined.
By Theorem~1 \cite{GR}, $\FF$ is an invariant hypersurface of the vector field
\begin{equation}
X \biggl( \frac{\pa}{\pa x} +  p \frac{\pa}{\pa y} \biggr) +
\sum_{i=2}^{m} Y_i \frac{\pa}{\pa p_i},
\label{91}
\end{equation}
which is obtained from \eqref{9} by eliminating the common factor $F^{-\mu}$.
By Lemma \ref{L1}, integral curves of \eqref{91} coincide with integral curves of \eqref{8},
hence $X=H$, $Y_i = H_i - p_iH_1$, and $\FF$ is an invariant hypersurface of \eqref{8}.
We remark that every solution of \eqref{6} is obtained from a solution of \eqref{8}
by choosing an appropriate parametrization $x_1(t)$. Thus $\FF$ is an invariant hypersurface
of equation \eqref{6} also.

Finally, consider the case $m=2$. Then $\dim PTM =3$ and $\dim \FF = 2$.
The field of contact planes $y_2=p_2 dx$ cuts a direction field on $\FF$, which coincides with
the restriction of the field \eqref{8} to $\FF$. Hence the projection of integral curves from
the surface $\FF \subset PTM$ to $M$ are isotropic lines
and non-parametrized extremals of the functional $J^{(n)}$ simultaneously.
\ebox

In accordance with what the previous reasonings, one can give the following definition.

\begin{defin}
\label{DEF1}
The projections of integral curves of equation \eqref{8} from $PTM$ to $M$
distinguished from a point are non-parametrized geodesics in the pseudo-Finsler space $(M, \oF)$.
\end{defin}

By Lemma~\ref{L2}, in the case $m=2$ all isotropic lines are geodesics in the sense of the given definition.

The natural parametrization of non-isotropic geodesics is defined by equation \eqref{5} with
$\psi(s)=s^{\frac{2}{n}}$ and coincides with the arc-length parametrization.
In the case $n>2$ the natural parametrization of isotropic geodesics is not defined, while
in the case $n=2$ it is defined by equation \eqref{5} with $\psi = \id$.

\section{Polynomial pseudo-Finsler metrics on 2-manifolds}\label{sect2}

From now on, we consider the case when $m=2$ and the function $\oF$ is a homogeneous polynomial
of degree $n \ge 2$ in $(\dx_i)$. Denote the coordinates on the manifold $M$ by $(x,y)$.

Consider pseudo-Finsler space with the metric function $\of = \oF^{\frac{1}{n}}$, where
\begin{equation}
\oF(x,y;\dx,\dy) = \sum_{i=0}^{n} a_i(x,y) \dx^{n-i}\dy^i, \ \ \
F(x,y;p) = \sum_{i=0}^{n} a_i(x,y) p^i,
\label{10}
\end{equation}
the coefficients $a_i$ smoothly depend on $(x,y)$.
Then equation \eqref{7} reads
\begin{equation}
\ddx = \frac{\oH_1}{\oH}, \ \ \ddy = \frac{\oH_2}{\oH}, \ \ \
\oH =
\left|\begin{matrix}
\oF_{\dx\dx} & \! \oF_{\dx\dy} \phantom{\bigl|}\!\! \\
\oF_{\dx\dy} & \! \oF_{\dy\dy} \phantom{\bigl|}\!\! \\
\end{matrix}\right|,
\ \,
\oH_1 =
\left|\begin{matrix}
\oG_1 & \! \oF_{\dx\dy} \phantom{\bigl|}\!\! \\
\oG_2 & \! \oF_{\dy\dy} \phantom{\bigl|}\!\! \\
\end{matrix}\right|,
\ \,
\oH_2 =
\left|\begin{matrix}
\oF_{\dx\dx} & \! \oG_1 \phantom{\bigl|}\!\! \\
\oF_{\dx\dy} & \! \oG_2 \phantom{\bigl|}\!\! \\
\end{matrix}\right|,
\label{11}
\end{equation}
where $\oG_1 = \oF_x - \dx \oF_{\dx x} - \dy \oF_{\dx y}$ and $\oG_2 = \oF_y - \dx \oF_{x \dy} - \dy \oF_{\dy y}$.

\begin{lemma}
\label{L3}
The projectivization $\Pi \colon T\M \to PT\M$ sends equation \eqref{11} to
\begin{equation}
p = \frac{dy}{dx}, \ \ \ \frac{dp}{dx} = \frac{H_2 - pH_1}{H} = \frac{P}{\Del},
\label{12}
\end{equation}
where
\begin{equation}
\label{13}
\begin{aligned}
&\Del (x,y;p) = nFF_{pp} - (n-1) F^2_p, \\
& P (x,y;p) = nF (F_y-F_{xp}-pF_{yp}) + (n-1)F_p (F_x+pF_y). \\
\end{aligned}
\end{equation}
\end{lemma}

{\bf Proof}.
Taking into account \eqref{8}, it remains to establish the equality
$\frac{H_2 - pH_1}{H} = \frac{P}{\Del}$, where $\Del, P$ are defined in \eqref{13}.
Let us prove that $H = (n-1)\Del$ and $H_2 - pH_1 = (n-1)P$, i.e.,
$\oH = \dx^{2n-4}(n-1)\Del$ and $\oH_2 - p\oH_1 = \dx^{2n-2}(n-1)P$.
Since both sides of the two last equalities can be treated as quadratic forms on
$a_0, \ldots, a_n$ with coefficients depending on $\dx, \dy$, it suffices to compare
the coefficients of the monomials $a_ia_j$ in the lift and right-hand sides.

Put $\eps_{ij}=1$ if $i \neq j$ and  $\eps_{ij}=\frac{1}{2}$ if $i = j$.
Direct calculation shows that the coefficient of the monomial $a_ia_j$, $i+j=k$, in the expression
$\oH = \oF_{\dx\dx}\oF_{\dy\dy}- \oF_{\dx\dx}^2$ is $\alf_{ij}\eps_{ij} \dx^{2(n-1)-k}\dy^{k-2}$, where
\begin{multline}
\alf_{ij} =
(n-i)(n-1-i)j(j-1) + (n-j)(n-1-j)i(i-1) - 2ij(n-i)(n-j) = \\
(n-1) \bigl(n(k^2-k-4ij)+2ij \bigr).
\label{14}
\end{multline}
On the other hand, the coefficient of the monomial $a_ia_j$, $i+j=k$, in the expression
$\Del = nFF_{pp} - (n-1) F^2_p$ is $\bet_{ij}\eps_{ij} p^{k-2}$, where
\begin{equation}
\bet_{ij} = n (i(i-1) + j(j-1)) - 2ij (n-1) = n(k^2-k-4ij)+2ij.
\label{15}
\end{equation}
From \eqref{14} and \eqref{15}, we have $\alf_{ij}=(n-1)\bet_{ij}$,
that proves $\oH = \dx^{2n-4}(n-1)\Del$.
The proof of the equality $\oH_2 - p\oH_1 = \dx^{2n-2}(n-1)P$ is similar.
\ebox

\begin{remark}
\label{R1}
{\rm
From formula \eqref{13} it follows that $\Del$ and $P$ are polynomials in $p$
of degrees not greater than $2n-4$ and $2n-1$, respectively. For instance,
$$
\Del(x,y;p) = (2na_na_{n-2} - (n-1)a_{n-1}^2) p^{2n-4} + \cdots + 2na_0a_2 - (n-1)a_1^2.
$$
}
\end{remark}


For our further purposes, it is convenient to write equation \eqref{12} as the field
\begin{equation}
\Del \biggl( \frac{\pa}{\pa x} +  p \frac{\pa}{\pa y} \biggr) + P \frac{\pa}{\pa p}.
\label{16}
\end{equation}
The field \eqref{16} is defined in the whole space $PTM$ including the isotropic surface $\FF$.
The field of contact planes $dy=pdx$ defines on $\FF$ a direction field whose integral curves correspond to
isotropic lines, while all remaining integral curves of the field \eqref{16}
(that do not belong entirely to the isotropic surface)
correspond to non-isotropic geodesics.

In accordance with Definition~\ref{DEF1}, non-parametrized geodesics in the pseudo-Finsler space $(M, \oF)$ are
the projections of integral curves of the field \eqref{16} from $PTM$ to $M$
distinguished from a point.
Singularities of geodesics occur at the points of $PTM$ where $\Del(x,y;p)$ vanishes.
To describe the locus of such points, we use the following lemma.

\begin{lemma}
\label{L4}
Given a polynomial
\begin{equation}
\F(p)= \prod_{i=1}^n (p+\gam_i),  \ \ \, \gam_i \in \bR, \ \ n \ge 2,
\label{17}
\end{equation}
consider the polynomial
\begin{equation}
\Del (p) = n \F(p) \F''(p) - (n-1) \F'(p)^2.
\label{18}
\end{equation}
Then the following statements hold:
\begin{itemize}
\item[{\rm (a)}]
$\Del \equiv 0$ if and only if $\gam_1 = \cdots = \gam_n$.
\item[{\rm (b)}]
Suppose that $\gam_i \neq \gam_j$ for at least one pair $i,j$.
Then $p$ is a real root of the polynomial $\Del$ if and only if $p$ is a multiple root of the polynomial $\F$.
\item[{\rm (c)}]
If $p$ is a double root of $\F$ and $n \ge 3$, then $p$ is a double root of $\Del$.
\end{itemize}
\end{lemma}

{\bf Proof}.
The implications
$\gam_1 = \cdots = \gam_n$ $\Rightarrow$ $\Del \equiv 0$ in (a) and
$\F(p) = \F'(p) = 0$ $\Rightarrow$ $\Del (p) = 0$ in (b) are trivial.
The implication $\Del \equiv 0$ $\Rightarrow$ $\gam_1 = \cdots = \gam_n$ in (a)
follows from (b). Indeed, assume that $\Del \equiv 0$ holds and
any two of the numbers $\gam_1, \ldots, \gam_n$ are not equal.
By~(b), $\Del (p) = 0$ implies $\F(p) = 0$. Hence $\F \equiv 0$, which contradicts \eqref{17}.

Statement (c) is also trivial: differentiating \eqref{18} twice, from
$\F(p)=\F'(p)=0$ we get $\Del(p)=\Del'(p)=0$ and $\Del''(p)=(2-n){\F''(p)}^2 \neq 0$
if $\F''(p) \neq 0$ and $n \neq 2$.

It remains to prove the implication
$\Del (p) = 0$ and $p \in \bR$ $\Rightarrow$ $\F (p) = \F' (p) = 0$ in the statement (b).
Assume that $\gam_i \neq \gam_j$ for at least one pair $i,j$ and
there exist $p_* \in \bR$ such that $\F (p_*) \neq 0$.
Making the change of variables $p \mapsto p-p_*$, without loss of generality we can assume that $p_*=0$.
Then $\F(0) = \gam_1 \cdots \gam_n \neq 0$, and
$$
\F'(0) = \sum_{i=1}^n \alf_i \F(0), \ \ \
\F''(0) = 2\sum_{i<j} \alf_i \alf_j \F(0), \ \ \ \alf_i = \frac{1}{\gam_i}.
$$
Substituting the above formulae in \eqref{18}, after straightforward transformations we get
\begin{equation}
\Del (0) = \F^2(0) \Bigl( 2n \sum_{i<j} \alf_i \alf_j - (n-1) \Bigl( \sum_{i=1}^n \alf_i \Bigr)^2 \Bigr) =
-\F^2(0) \phi_n (\alf_1, \ldots, \alf_n),
\label{24-08-2016}
\end{equation}
where
\begin{equation*}
\phi_n (\alf_1, \ldots, \alf_n) = n \sum_{i=1}^n \alf_i^2 - \Bigl( \sum_{i=1}^n \alf_i \Bigr)^2.
\end{equation*}

Let us prove that for any $n \ge 2$ the form $\phi_n (\alf_1, \ldots, \alf_n) \ge 0$
and $\phi_n (\alf_1, \ldots, \alf_n)=0$ if and only if $\alf_1 = \cdots = \alf_n$.
Indeed, consider the vectors $\alf=(\alf_1, \ldots, \alf_n)$ and $\bet = (1, \ldots, 1)$
in $n$-dimensional Euclidean space with the standard inner product.
Then the Cauchy--Schwarz inequality $(\alf,\bet)^2 \le (\alf,\alf) (\bet,\bet)$
gives the required assertion.

By our assumption $\alf_i \neq \alf_j$ for at least one pair $i,j$.
Then $\phi_n (\alf_1, \ldots, \alf_n) >0$, and equality \eqref{24-08-2016}
implies that $\Del (0) = 0$  $\Leftrightarrow$ $\F (0) = 0$.
Moreover, from \eqref{18} it follows $\Del (p_*) = \F (p_*) = 0$ $\Rightarrow$ $\F' (p_*) = 0$, i.e.,
$p_*=0$ is a multiple root of $\F$. The lemma is proved.
\hfill\ebox

\begin{remark}
\label{R2}
{\rm
Obviously, the implication $\F = \F' = 0$ $\Rightarrow$ $\Del = 0$ holds true for all polynomial $\F$,
not necessarily \eqref{17}. However, the inverse implication $\Del = 0$ $\Rightarrow$ $\F = 0$
is not valid if $\F$ has a complex root.
The reason for this is easily ascertained: the inequality $\phi_{n} (\alf_1, \ldots, \alf_n) \ge 0$ is not valid
if among the numbers $\alf_i$ some are complex.

For example, consider the polynomial $\F = p^3+p$ with the unique real root $p=0$.
Then the corresponding polynomial $\Del = 2(3p^2-1)$ has two real roots,
none of those coincides with $0$.
Moreover, the polynomial $\F = p^4+6p^2+1$ does not have real roots at all, while
the corresponding polynomial $\Del = 48(p^2-1)^2$ has two double roots $p= \pm 1$.
}
\end{remark}

Lemma~\ref{L4} gives a simple geometrical description of the singular locus of equation \eqref{12}
for the domain $M' \subset M$,
where pseudo-Finsler space $(M, \oF)$ has $m$ isotropic lines passing trough every point of $M'$, i.e.,
the polynomial $F(p)$ has $m$ real roots (taking into account the multiplicity and possibly including $p=\infty$).
For $(x,y) \in M'$ the function $\Del (x,y;p)$ vanishes if and only if
at least two of $m$ isotropic lines are tangent at $(x,y)$
and $p$ is the corresponding tangential direction.
Remark that this statement is not valid for the complement of $M'$, where the polynomial $F(p)$ has complex roots.

This question will be considered in detail for $n=3$.

\subsection{Pseudo-Riemannian metrics}\label{sect2.1}

By Remark \ref{R1}, in the case $n=2$ (pseudo-Riemannian metrics) $\Del$ is a zero degree polynomial in $p$,
that is, $\Del$ does not depend on $p$. Moreover, it is easy to check that $\Del = -4 D_{[F]}$,
where $D_{[F]}$ means the discriminant of the quadratic polynomial $F$. Hence the locus of singularities
of equation \eqref{12} coincides with the discriminant curve of the implicit differential equation $F(x,y;p)=0$.
It is not hard to see that the equation $\Del(x,y)=0$ defines an invariant surface of the field \eqref{16}
filled with integral curves whose projections $PTM \to M$ are points forming the discriminant curve.

This property leads to a curious phenomenon: geodesics cannot pass through a point $(x,y)$ of the discriminant curve
in arbitrary tangential directions, but only in admissible directions $p$ defined by the condition $P(x,y;p)=0$.
Generically, $P(x,y;p)$ is a cubic polynomial in $p$ and the number of admissible directions is 1 or 3.
Singularities of the geodesic flows in pseudo-Riemannian metrics are studied in detail
(the interested reader is referred to the papers \cite{GR, Rem-Pseudo, Rem15} devoted to 2-dimensional
pseudo-Riemannian metrics; similar results for 3-dimensional pseudo-Riemannian metrics were announced in \cite{PR}).

It should be remarked that the case $n=2$ is exceptional from the viewpoint of Finsler and pseudo-Finsler
geometry ($n>2$).
In the case $n>2$, $\Del$ generically depends on $p$ and the notion of admissible directions does not appear.
Geodesics pass through every point of $M$ in all possible directions,
but some directions at some points are singular.
In other words, only points of the space $PTM$ may have the property of being singular.

In the rest of the paper, we consider the case $n=3$ (cubic pseudo-Finsler metrics) in detail.

\subsection{Cubic pseudo-Finsler metrics}\label{sect2.2}

Let $D_{[F]}$ and $D_{[\Del]}$ be the discriminants of the cubic polynomial $F(x,y;p)$
and the quadratic polynomial $\Del(x,y;p)$ in $p$, respectively. A direct calculation shows that
$D_{[\Del]} = -12D_{[F]}$.

Introduce the following stratification of the manifold $M$.
The open domains $M_+, M_-$ are defined by the conditions $D_{[F]}>0$, $D_{[F]}<0$, respectively.
Generically, $M_+, M_-$ are separated by the discriminant curve $M_0 \colon D_{[F]}=0$,
which consists of regular points (the cubic polynomial $F$ has one prime root and one double root)
and cusps ($F$ has a triple root). By $M_{0,1}$ denote the set of all regular points of $M_0$,
while $M_{0,0} = M_{0} \setminus M_{0,1}$.
The discriminant of the quadratic polynomial $\Del$ is strictly negative in $M_+$, hence
singularities of equation \eqref{12} occur only in $M_-$ and $M_{0}$.
Further we exclude from consideration the stratum $M_{0,0}$ of dimension zero,
and consider only $M_-$ and $M_{0,1}$.

In a neighborhood of every point of $M \setminus M_{0,0}$
the cubic polynomial $F$ has at least one prime real root $p_*(x,y)$ smoothly depending on $x,y$.
To simplify calculations, choose local coordinates such that the integral curves
of the vector field $\frac{dy}{dx}=p_*(x,y)$ (one of three families of isotropic lines)
become $x=\const$. This yields $a_3(x,y) \equiv 0$ and
\begin{equation}
\label{19}
\begin{aligned}
& F = ap^2 + 2bp + c, \ \ \,  \Del = -2 (ap+b)^2 + 6(ac-b^2), \ \ \,  D_{[F]} = 4a^2 (b^2-ac),
\\
& M_{\pm} =\{\pm (b^2-ac)>0, \ a \neq 0\}, \ \ \, M_{0,1} =\{b^2-ac=0, \ a \neq 0\},
\\
& P = 3F (F_y-F_{xp}-pF_{yp}) + 2F_p (F_x+pF_y).
\end{aligned}
\end{equation}

\subsubsection{Singularities in the stratum $M_-$}\label{sect2.2.1}

At every point in $M_-$, the quadratic equation $\Del =0$ has two prime real roots
\begin{equation}
p_{1,2} = \frac{\pm \del - b}{a}, \ \ \, \del = \sqrt{3(ac-b^2)},
\label{20}
\end{equation}
and the domain $M_-$ is filled with two transverse families of integral curves
of the binary implicit differential equation $\Del=0$,
which we shall call {\it singular lines} of the metric.

Consider the curves $S_i \subset M_-$ defined by the equations $P(x,y;p_i) = 0$, $i=1,2$,
where $P$ is defined in \eqref{19}.
They can be also considered as the branches of the locus $\operatorname{res} (\Del,P)=0$,
where ``$\operatorname{res}$'' means the resultant of two polynomials in $p$.
In the space $PTM$, consider the corresponding curves
$$
\ov S_i = \{(x,y;p_i) : (x,y) \in S_i\}, \ \ i=1,2,
$$
which consist of singular points of the field \eqref{16}.

By $\Gam_q$ denote the family of geodesics outgoing from a point $q=(x,y)$.
The simplest type of singularities of the geodesic flow (codimension 0) is given in the following theorem.

\begin{theorem}
\label{T1}
Let $q \in M_-$ and $(q;p_i) \notin \ov S_i$, $i=1,2$.
Then there exists a unique geodesic passing through the point $q$ with tangential direction $p_i$:
a semicubic parabola with the cusp at $q$.
In particular, if $q \in M_- \setminus (S_1 \cup S_2)$, the family $\Gam_q$ contains two semicubic parabolas
with tangential directions $p_i$,
while geodesics with all remaining directions at $q$ are smooth.
\end{theorem}

{\bf Proof}.
If $P(q;p_i) \neq 0$, then, by the standard existence and uniqueness theorem,
the field \eqref{16} has a unique integral curve $\gam_i$ passing through the point $(q;p_i)$.
From the conditions $a \neq 0$, $ac-b^2 \neq 0$ it follows that $\Del$ and $\Del_p$ do not vanish simultaneously,
see \eqref{19}.

Hence the curve $\gam_i$ has the first order tangency with the {\it vertical} direction
(the {vertical} direction in the space $PTM$ is called the $p$-direction, i.e.,
the kernel of the natural projection $PTM \to M$)
at the point $(q;p_i)$, and the projection of the curve $\gam_i$
to $M$ is a semicubic parabola with the cusp at $q$.
\ebox

\begin{figure}[ht]
\begin{center}
\includegraphics[height=5.4cm]{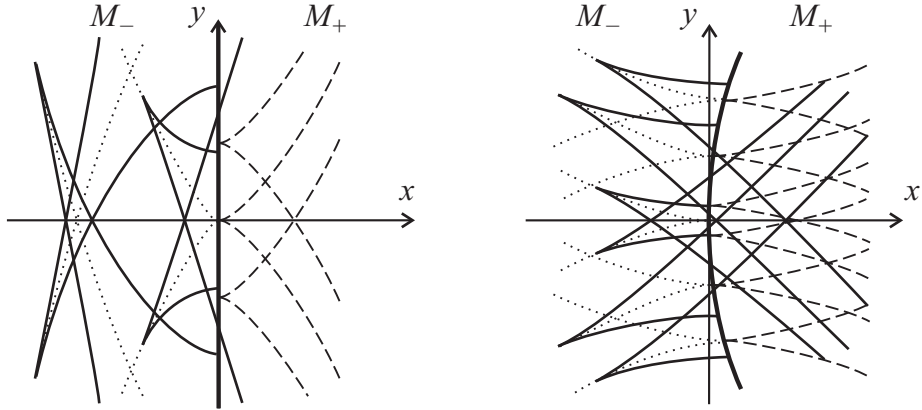}
\end{center}
\caption{
Example \ref{Example1}: I and II (left and right, resp.).
The stratum $M_{0,1}$ (the $y$-axis on the left and the parabola on the right) is depicted with the bold solid line.
Solid and dotted curves present geodesics and singular lines, respectively.
Dashed curves present isotropic lines.
}
\label{fig1}
\end{figure}

\begin{exam}
{\rm
Let $F$ be given by formula \eqref{19} with $a=1$ and $b=0$.
Then $\Del = 2(3c-p^2)$ and $P = 7c_y p^2 + 4c_xp + 3cc_y$.
The stratums $M_-$ and $M_{0,1}$ are defined by the conditions $c(x,y)>0$ and $c(x,y)=0$, respectively,
and the curves $S_i$ are $c_x = \pm 2c_y \sqrt{3c}$.

{\bf I.} \
Put $c(x,y)=-x$ (Fig.~\ref{fig1}, left).
Then $S_1=S_2=\emptyset$ and the semiplane $x > 0$ ($M_+$) is filled with the net of isotropic lines
$y = \pm \frac{2}{3}x^{\frac{3}{2}} + \const$ (dashed curves), while
the semiplane $x <0$ ($M_-$) is filled with the net of singular lines
$y = \pm \frac{2}{\sqrt{3}}(-x)^{\frac{3}{2}} + \const$
(dotted curves).
Cusps appear when geodesics (solid curves) are tangent to singular lines.
Remark that geodesic pass from $M_-$ to $M_+$ or vise versa through $M_{0,1}$ (the $y$-axis)
without singularity if they intersect the $y$-axis with any non-isotropic tangential direction $p \neq 0$.
Otherwise, equation \eqref{12} has singularity. As we shall see in Section~\ref{sect2.2.2}, at such points
there exist a one-parameter family of geodesics outgoing in both domains $M_+$ and $M_-$ with the common
tangential direction $p=0$, and the prolongation of geodesics through $M_{0,1}$ is not naturally defined.

{\bf II.} \
Put $c(x,y)=\alf y^2-x$ with $\alf \neq 0$.
Then $S_i \neq \emptyset$, $i=1,2$, but both curves $S_i$ do not pass through the origin.
In a neighborhood of the origin that does not contain the curves $S_i$, geodesics
are presented in Fig.~\ref{fig1} (right).
Here, for definiteness, we assume $\alf>0$. All notations have the same meanings as before.
}
\label{Example1}
\end{exam}

The next type of singularities of the geodesic flow in the domain $M_-$ (codimension 1) is connected with
vanishing of the field \eqref{16}.
This field belongs to a special class of vector fields
whose singular points are not isolated, but form a manifold $W$ of codimension~2 in the phase space.
Yet such fields appear in many problems, see e.g. \cite{AGiv, GR, Martinet, Rem06, Rem08, Rouss, ZhitSo}.
It is convenient to expressed the above condition in the following algebraic form:
the germs of all components of the field at every singular point belong to the ideal $I$
(in the ring of smooth germs) generated by two of them.

The spectrum of the linear part of such a field (for brevity, we shall call it {the spectrum of the field})
contains only two non-zero eigenvalues $\la_{1,2}$, which play a prominent role in establishing the local
normal form of the field.
For instance, all components of the field \eqref{16} belong to the ideal $I = \<\Del,P\>$,
and the set of singular points $W = \ov S_1 \cup \ov S_2$.
The eigenvalues $\la_{1,2}$ and the corresponding eigenvectors are described by the following lemma.

\newpage

\begin{lemma}
\label{L5}
$\phantom{.}$

1.
The resonance $\la_1+\la_2 = 0$ holds at all points $(q;p) \in \ov S_i$, i.e.,
$\la_{1,2}$ are real or pure imaginary numbers with opposite signs.

2. The following conditions are equivalent:

2.1.
The eigenvalues $\la_{1}, \la_{2}$ at $(q;p) \in \ov S_i$ are not equal to zero.

2.2.
$\ov S_i$ is a regular curve transversal to the contact plane $pdx - dy = 0$ at $(q;p)$.

2.3.
$S_i$ is a regular curve and the direction $p_i$ is transversal to $S_i$ at $q$.

3.
Generically, at almost all points $(q;p) \in \ov S_i$ the conditions 2.1\,--\,2.3 hold.
\end{lemma}

{\bf Proof}.
Without loss of generality, suppose that $q=0$ (the origin) belongs to $S_1$
and choose local coordinates centered at $0$
that preserve the lines $x=\const$ and send integral curves of the vector field $\frac{dy}{dx}=p_1(x,y)$
to parallel lines $y=\const$.
The existence of such local coordinates follows from
the general fact: if $V_1$ and $V_2$ are smooth vector fields on the plane transversal at the point $0$, then
in a neighborhood of $0$ there exist local coordinates such that integral curves of $V_1$ and $V_2$ coincide
with the coordinate lines.

Then the polynomials $F, \Del, D_{[F]}$ have the form \eqref{19}
and the identities $p_1 \equiv 0$ and $b(x,y) \equiv \del(x,y)$ hold.
Note that the first of them implies $3ac \equiv 4b^2$.
From $b^2-ac<0$ it follows $ac > 0$. Taking into account $3ac \equiv 4b^2$,
we conclude that none of the coefficients $a,b,c$ vanishes at $0$.
Below, we present the proof for the stratum $S_1$ given by the equation $3c(c_y-2b_x) +4bc_x = 0$.
The proof for the stratum $S_2$ is similar.

Let $\Lambda$ be the matrix of the linear part of the field \eqref{16}
and  $\Lambda_1$ be the matrix of the Pfaffian system $d \Del =0$, $d P = 0$, $pdx - dy = 0$
considered at arbitrary point $(q;p) \in \ov S_1$, that is, for $q \in S_1$ and $p=0$:
\begin{equation*}
\Lambda =
\left(\begin{matrix}
\Del_x  &  \Del_y  &  \Del_p \\
0  &  0  &  0 \\
P_x  &  P_y  &  P_p \\
\end{matrix}\right)
\IIII_{p=0},
\qquad
\Lambda_1 =
\left(\begin{matrix}
\Del_x  &  \Del_y  &  \Del_p \\
P_x  &  P_y  &  P_p \\
0  &  -1  &  0 \\
\end{matrix}\right)
\IIII_{p=0},
\end{equation*}
where
\begin{equation}
\begin{aligned}
& \Del_x \I_{p=0} = 6(ac_x+a_xc) - 16bb_x, \ \ \
   \Del_p \I_{p=0} = -4ab, \\
& P_x \I_{p=0} = c_x(3c_y-2b_x) + 4bc_{xx} +  3c (c_{xy}  - 2b_{xx}), \\
& P_p \I_{p=0} = 4ac_x - 6a_xc  + 2b (5c_y - 2b_x). \\
\end{aligned}
\label{21}
\end{equation}

1. \
To prove the first statement, it suffices to show that $\tr \Lambda = 0$.
Taking into account the equality $3c(c_y-2b_x) +4bc_x = 0$ on $S_1$ and
the identity $b \equiv \del$ (that implies $3ac \equiv 4b^2$ in a neighborhood of the origin),
from \eqref{21} we have
$$
\tr \Lambda =
(\Del_x+P_p) \I_{p=0} = 10 [ac_x + b(c_y-2b_x)] =  10 \Bigl( ac_x - b\frac{4bc_x}{3c}\Bigr) =
\frac{10c_x}{3c} (3ac-4b^2) = 0.
$$
Hence the characteristic equation of the matrix $\Lambda$ reads $\la (\la^2 + |\Lambda_1|) = 0$,
this yields the equation  $\la^2 + |\Lambda_1| = 0$ for $\la_{1,2}$.

2. \
Differentiating the identity $4b^2 \equiv  3ac$ by $x$, we get $8bb_x \equiv 3(a_xc+ac_x)$.
Using these identities and \eqref{21}, we have
\begin{multline*}
\Del_x \I_{p=0} = 6(ac_x+a_xc) - 16 bb_x = 6(ac_x+a_xc) - 6(ac_x+a_xc) = 0, \\
|\Lambda_1| = - \Del_p P_x \I_{p=0} =
16ab^2c_{xx} - 8abb_xc_x + 12ab [c_xc_y + cc_{xy} - 2b_{xx}c] =\\
12a^2cc_{xx} - 3ac_x(a_xc+ac_x) + 12ab [c_xc_y + cc_{xy} - 2b_{xx}c].
\end{multline*}
Thus the condition $\la_{1,2} \neq 0$ is equivalent to $|\Lambda_1| \neq 0$ that, in turn,
is equivalent to the condition~2.2.

On the other hand, the curve $S_1$
is tangent to the direction $p=0$ at the point $q=0$ if and only if $[3c(c_y-2b_x) +4bc_x]'_x = 0$.
Taking into account the equalities $4b^2 \equiv  3ac$ and $8bb_x \equiv 3(a_xc+ac_x)$, we get
\begin{multline*}
[3c(c_y-2b_x) +4bc_x]'_x =
3 (c_xc_y + cc_{xy}) - 6b_{xx}c + \frac{4b^2c_{xx} - 2bb_xc_x}{b} = \\
3 (c_xc_y + cc_{xy}) - 6b_{xx}c + \frac{12acc_{xx} - 3(a_xc+ac_x)c_x}{4b} =
\frac{|\Lambda_1|}{4ab}.
\end{multline*}
This proves that $\la_{1,2} \neq 0$ is equivalent to the condition~2.3.

3. \
Generically, at almost all points $(q;p) \in \ov S_1$ the determinants
\begin{equation*}
\left|\begin{matrix}
\Del_x  &  \Del_y  \\
  P_x   &    P_y   \\
\end{matrix}\right|,
\quad \
|\Lambda_1| =
\left|\begin{matrix}
\Del_x  &  \Del_p  \\
  P_x   &    P_p   \\
\end{matrix}\right|
\end{equation*}
are not equal to zero. Hence $S_1$ and $\ov S_1$ are regular curves and, moreover,
the conditions 2.1\,--\,2.3 hold.
\ebox

\begin{theorem}
\label{T2}
Let $(q;p_i) \in \ov S_i$ be a generic singular point of the field \eqref{16}.
Then the germ \eqref{16} at $(q;p_i)$ is smoothly orbitally equivalent to
\begin{eqnarray}
&&
\xi \frac{\pa}{\pa \xi} - \eta \frac{\pa}{\pa \eta} + \xi\eta \frac{\pa}{\pa \zeta},
\ \ \ \textrm{if} \ \ \, \la_{1,2} \in \bR \setminus 0,
\label{22}
\\
&&
\eta \frac{\pa}{\pa \xi} - \xi \frac{\pa}{\pa \eta} + (\xi^2+\eta^2) \frac{\pa}{\pa \zeta},
\ \ \ \textrm{if} \ \ \, \la_{1,2} \in \bI \setminus 0,
\label{23}
\end{eqnarray}
where $\bR, \bI$ are real and imaginary axes, respectively.

In the first case there exist two geodesics passing through the point $q \in S_i$ with the tangential direction $p_i$,
both of them smooth. In the second case, there are no geodesics passing through the point $q \in S_i$
with the tangential direction $p_i$.
\end{theorem}

{\bf Proof}.
Since $(q;p_i) \in \ov S_i$ is a generic singular point, $|\Lambda_1| \neq 0$.
By Lemma~\ref{L5}, the eigenvalues $\la_{1,2} \neq 0$ and $\ov S_i$ is a regular curve consisting
of singular points of the field \eqref{16}.
The linear part of the germ \eqref{16} at $(q;p_i)$
(and every singular point sufficiently close to $(q;p_i)$)
is orbitally equivalent to the linear part of the field \eqref{22} or \eqref{23} if
$|\Lambda_1| < 0$ or $|\Lambda_1| > 0$, respectively.
Here we use the following terminology:
two vector fields are called orbitally smoothly (resp. topologically) equivalent,
if there exists a diffeomorphism (resp. homeomorphism) that conjugates their integral curves, i.e.,
orbits of their phase flows.\footnote{
It slightly differs from the generally accepted definition of the orbital equivalence,
where coincidence of the orientation of integral curves is also required.
Our definition is naturally related to directions fields, whose
integral curves do not have an orientation a priori.
}

Recall that the field \eqref{16} belongs to the class of vector fields
whose singular points are not isolated, but form a manifold $W$ of codimension~2 in the phase space
(in our case, $W=\ov S_i$). Local normal forms of such fields were studied by many authors.
In \cite{RemizovTari} (Appendix A), we present a brief survey of such results,
which covers all cases with $\Real \la_{1,2} \neq 0$.
This condition is equivalent to the assumption that $W=\ov S_i$ is the local center manifold,
and consequently, the phase portrait of \eqref{16} has a simple topological structure
(we shall discuss it later on, in the proof of Theorem~\ref{T4}).
For instance, in the case $\la_{1,2} \in \bR \setminus 0$, the germ \eqref{16} with generic quadratic part
is smoothly orbitally equivalent to \eqref{22}. This result belongs to Roussarie \cite{Rouss}.
The genericity condition is determined explicitly in \cite{RemizovTari} (Theorem 5.7).
The case $\la_{1,2} \in \bI \setminus 0$ is more complicated.
However, in \cite{AGiv} (Chapter~2, Section~1.2) and \cite{Martinet} it claims that
in this case the germ \eqref{16} with generic quadratic part is smoothly orbitally equivalent to \eqref{23}.

Remark that the diffeomorphism that brings the germ \eqref{16} to the normal form \eqref{22} or \eqref{23}
does not give a normal form of equation \eqref{12}, since it does not preserve the contact structure $dy = pdx$.
However, we need not a normal form of \eqref{12}.

To prove the last statement of the theorem, we only need to consider the possible mutual relationships between
the phase portrait of the germ \eqref{16} and the $(x,y)$-plane in the space $PTM$.
Geodesics are obtained from those integral curves of the field \eqref{16} whose projection on the $(x,y)$-plane
are distinguished from points. Moreover, isotropic geodesics correspond to those integral curves
that belong to the isotropic hypersurface $\FF$
(by Lemma~\ref{L2}, $\FF$ is an invariant hypersurface of the field \eqref{16}).

We consider the real and imaginary cases separately.

{\it The real case}.
The field \eqref{22} has the first integral $\xi\eta$.
The invariant foliation $\xi\eta = \const$ contains only two leaves $\xi=0$ and $\eta=0$ that
pass through singular points of the field, while all remaining invariant leaves are hyperbolic cylinders
$\xi\eta = \const \neq 0$, which do not intersect the set of singular points.
It is easy to see that for every singular point of the field \eqref{22} there are only two integral curves
passing through this point: the straight lines parallel to the $\xi$-axis and the $\eta$-axis, respectively.

We prove now that the eigenvectors with the eigenvalues $\la_{1,2} \neq 0$ are not vertical.
Let $e$ be an eigenvector of the matrix $\Lambda$ with $\la_i$.
Then $\Lambda e = \la_i e$, and $e= \alf\pa_x + \bet\pa_p$, where
$(\Del_x\I_{p=0}-\la_i)\alf + \Del_p\I_{p=0}\bet = 0$. If the eigenvector $e$ is vertical,
i.e., $\alf=0$, $\bet \neq 0$, this equality yields $\Del_p\I_{p=0}=0$.
From \eqref{21}, we have $a(0)=0$ or $b(0)=0$.
This contradicts to the fact (established in the proof of Lemma~\ref{L5}) that
none of the coefficients $a,b,c$ vanishes at $0$.

\begin{figure}[ht]
\begin{center}
\includegraphics[height=6.0cm]{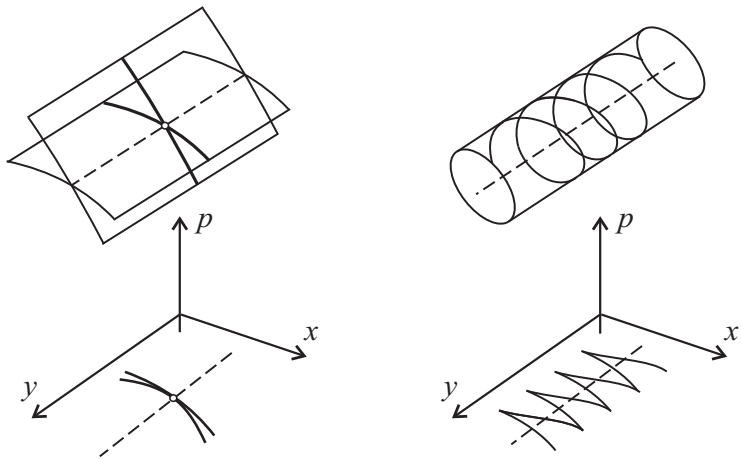}
\end{center}
\caption{
The phase portraits of the field \eqref{16} with the normal form \eqref{22} or \eqref{23} (left and right, resp.)
and the projections of its integral curves to $M$.
Dashed curves present $\ov S_i$ (up) and $S_i$ (down).
}
\label{fig2}
\end{figure}

From the considerations above, it follows that the field \eqref{16} has only two integral curves
passing through the given point $(q;p_i)$, both of them smooth and have non-vertical tangential directions.
Projecting these integral curves from $PTM$ to $M$, we get two smooth geodesics passing through
the point $q \in S_i$ with the tangential direction $p_i$; see Fig.~\ref{fig2} (left).

{\it The imaginary case}.
The field \eqref{23} has the first integral $\xi^2+\eta^2$.
The invariant foliation $\xi^2+\eta^2 = \const$ contains a one-dimensional degenerate leaf
$\xi=\eta=0$, which consists of singular points of the field \eqref{23}
and one-parameter family of two-dimensional leaves (elliptic cylinders $\xi^2+\eta^2 = \const \neq 0$),
which do not intersect the set of singular points.
The elliptic cylinders are filled with helix-like integral curves, whose
projections to $M$ have cusps; see Fig.~\ref{fig2} (right).

To complete the proof, observe that in both real and imaginary cases the curve $S_i$ itself is not a geodesic,
since $\ov S_i$ is transversal to the contact plane $pdx - dy = 0$ (statement~2.2 in Lemma~\ref{L5}).
Consequently, $\ov S_i$ is not a lift of a curve on $M$.
\ebox

\begin{figure}[ht]
\begin{center}
\includegraphics[height=4.8cm]{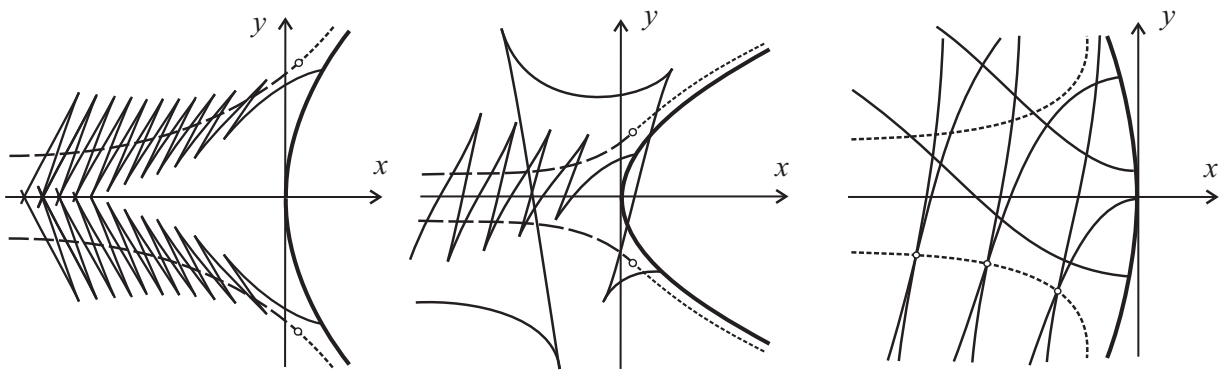}
\end{center}
\caption{
Example \ref{Example2}:
$F=p^2+c$, where $c=\alf y^2 - x$ with $\alf>0$ 
(left, center) and $\alf<0$ 
(right).
The solid lines present geodesics.
$M_{0,1}$ is depicted with bold solid line,
$S_i$ are depicted with short-dashed (long-dashed) lines
if $\la_{1,2} \in \bR \setminus 0$ ($\la_{1,2} \in \bI \setminus 0$, respectively).
}
\label{fig3}
\end{figure}

\begin{exam}
{\rm
Let $F$ be given by formula \eqref{19} with $a=1$, $b=0$, $c(x,y)=\alf y^2-x$.
Then $\Del = 2(3c-p^2)$, $P = 12\alf y p^2 - 4p + 6\alf yc$, and
the curves $S_i$ are the connected components of the graph $x = \alf y^2 - \frac{1}{48}/{\alf^2 y^2}$
lying the the upper and lower semiplanes.
A straightforward calculation shows that $p_i \I_{S_i} = 12\alf y (\alf y^2-x)$,
hence the direction $p_i$ is tangent to the curve $S_i$ at $x = \frac{47}{48}/{\sqrt{\alf}}$ only.
By the statement~2.3 in Lemma~\ref{L5}, the eigenvalues $\la_{1,2} \neq 0$ at all points of $\ov S_i$ if $\alf<0$
and at all points of $\ov S_i$ with $x \neq \frac{47}{48}/{\sqrt{\alf}}$ if $\alf \ge 0$.

In Fig.~\ref{fig3} (left and center) we present geodesics in the case $\alf>0$.   
Here both real and imaginary eigenvalues exist.
The parts of $S_i$ with real (imaginary) eigenvalues $\la_{1,2} \neq 0$
are presented as short-dashed (resp., long-dashed) lines.
The dots represent the points of $S_i$ with $x = \frac{47}{48}/{\sqrt{\alf}}$, where $\la_{1,2} = 0$.
In Fig.~\ref{fig3} (right) we present geodesics in the case $\alf<0$.   
Here only real eigenvalues exist, and the curves $S_i$ are presented as short-dashed lines.
The dots represent the points where geodesics intersect the curves $S_i$ with the singular tangential direction $p_i$.
}
\label{Example2}
\end{exam}

\subsubsection{Singularities in the stratum $M_{0,1}$}\label{sect2.2.2}

In this section as before,
we proceed in the local coordinates where $F$, $\Del$ and $D_{[F]}$ have the form \eqref{19}.
At every point $q \in M_{0,1}$ the polynomial $F$ has the double root $p_0 = -b/a$.
It is easy to see that $p_0$ is also a double root of the polynomial $\Del$ at $q$
(it follows from Lemma~\ref{L4} also).
Thus $(q; p_0)$, $q \in M_{0,1}$, are singular points of both implicit differential equations
$F=0$ and $\Del=0$. From \eqref{13} it follows that $P(q; p_0)=0$, hence $(q; p_0)$, $q \in M_{0,1}$,
are singular points of the field \eqref{16}.

Further we restrict ourselves to generic points $q \in M_{0,1}$ where
$M_{0,1}$ is a regular curve and the isotropic direction $p_0$ is transverse to $M_{0,1}$.
Then both implicit differential equations $F=0$ and $\Del=0$ have Cibrario normal form at such point
and their integral curves are semicubic parabolas lying on opposite sides of $M_{0,1}$
(in the domains $M_+$ and $M_-$, resp.) as it is presented in Fig.~\ref{fig1}.

\begin{figure}[ht]
\begin{center}
\includegraphics[height=5.4cm]{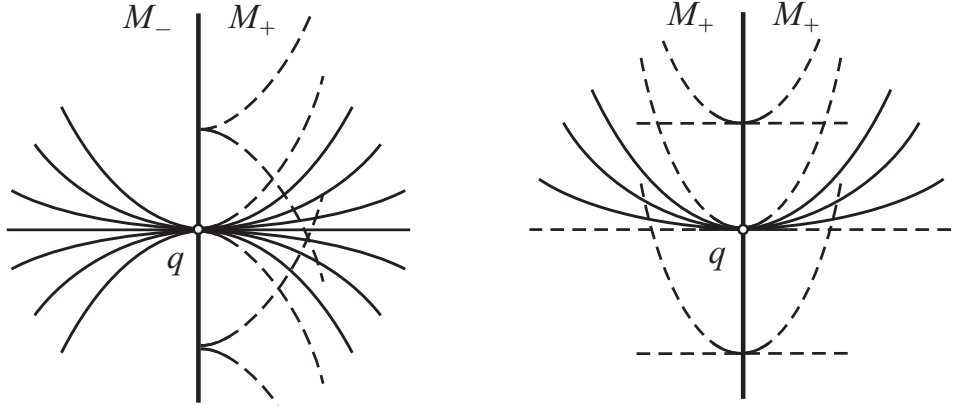}
\end{center}
\caption{
Illustrations to Theorem~\ref{T3} (left) and Theorem~\ref{T4} (right).
The stratum $M_{0,1}$ is depicted with the bold solid line.
Solid and dashed curves present non-isotropic and isotropic geodesics, respectively.
}
\label{fig4}
\end{figure}

\begin{theorem}
\label{T3}
Suppose that the isotropic direction $p_0$ is transverse to $M_{0,1}$ at $q \in M_{0,1}$.
Then the germ \eqref{16} at its singular point $(q; p_0)$ is smoothly orbitally equivalent to
\begin{equation}
3\xi \frac{\pa}{\pa \xi} + 2\eta \frac{\pa}{\pa \eta} + 0 \frac{\pa}{\pa \zeta},
\label{24}
\end{equation}
and to $p_0$ corresponds a one-parameter family of geodesics
outgoing from $q$ into $M_+$ and $M_-$.
There exist smooth local coordinates centered at $q$ such that this family is
\begin{equation}
x = \alf |\eta|^{\frac{3}{2}} + \eta^2 + \alf \bar X_{\alf}(\eta), \ \
y = \alf \eta |\eta|^{\frac{3}{2}} + \eps \eta^3 + \alf \bar Y_{\alf}(\eta), \ \ \eps \neq 0,
\label{25}
\end{equation}
where $\bar X_{\alf}(\eta) = o\bigl(|\eta|^{\frac{3}{2}}\bigr)$ and
$\bar Y_{\alf}(\eta) = o\bigl(|\eta|^{\frac{5}{2}}\bigr)$
are $C^2$-smooth functions.

Here $\alf>0$ $(\alf<0)$ corresponds to non-isotropic geodesics outgoing
from $q$ in $M_+$ (resp., $M_-$), while $\alf=0$ gives the isotropic geodesic,
a semicubic parabola lying in $M_+$. The limit case $\alf \to \infty$ corresponds to a unique
smooth geodesic passing through $q$ with the direction $p_0$.
In a neighborhood of $q$,
every non-isotropic geodesics outgoing from $q$ in $M_+$
belongs to the curvilinear tongue-like sector bounded by the branches of the isotropic geodesic
as it is presented in Figure~\ref{fig4} (left).
\end{theorem}

{\bf Proof}.
Without loss of generality, suppose that $q=0$ (the origin of the $(x,y)$-plane) and choose local
coordinates centered at $0$ that preserve the lines $x = \const$ and give $b(x,y) \equiv 0$.
It can be done using appropriate change of variables $y \mapsto yu(x,y)$, $u(0) \neq 0$.
Remark that (unlike Lemma~\ref{L5}) we do not have the identity $p_1 \equiv 0$ nor $p_2 \equiv 0$
in a neighborhood of $0$. Moreover, it is impossible to get any these identities using smooth change
of variables, since the integral curves of the implicit equation $\Del = 0$ with roots $p_1, p_2$
have cups on $M_{0,1}$.

In the local coordinates chosen above, we have
\begin{equation}
\begin{aligned}
&
F = ap^2 + c, \ \ \,
\Del = -2 (ap)^2 + 6ac, \ \ \,
D_{[F]} = -4a^3c,
\\
& P =
aa_y p^4 - 2aa_x p^3 + (7ac_y-3a_yc) p^2 + (4ac_x-6a_xc) p + 3cc_y.
\\
\end{aligned}
\label{26}
\end{equation}
The curve $M_{0,1}$ is given by $c(x,y)=0$ and the direction $p_0 = 0$ at every $q \in M_{0,1}$.
Hence the condition ``the direction $p_0$ is transverse to $M_{0,1}$ at $0$''
is equivalent to $c_x(0) \neq 0$.

Substituting $\Del$ and $P$ from \eqref{26} into \eqref{16}, one can find that the spectrum
of the field \eqref{16} at every point $(q; 0)$, $q \in M_{0,1}$, is $(\la_1,\la_2,0)$,
where $\la_1 = 6ac_x$, $\la_2 = 4ac_x$.
A straightforward computation shows that the corresponding eigenvectors are
\begin{equation}
e_1 = 2a\pa_x + 3c_y\pa_p, \ \
e_2 = \pa_p, \ \
e_0 = c_y\pa_x - c_x\pa_y.
\label{26-08-2016}
\end{equation}

Note that $\la_1:\la_2 \equiv 3:2$ at all points $(q; 0)$, $q \in M_{0,1}$,
the pair $(\la_1,\la_2)$ is non-resonant and belongs to the Poincar\'e domain.
Therefore, the germ \eqref{16} at $0$ is smoothly orbitally equivalent to the linear field \eqref{24}
(Theorem~5.5 in \cite{RemizovTari}).
Moreover, comparing \eqref{16} and \eqref{24}, one can see that the conjugating diffeomorphism
$(x,y,p) \mapsto (\xi,\eta,\zeta)$ can be chosen in the form
\begin{equation}
x = 2a\xi + c_y \zeta + f_1(\xi,\eta,\zeta), \ \
p = 3c_y\xi + \eta + f_2(\xi,\eta,\zeta), \ \
y = - c_x \zeta + f_3(\xi,\eta,\zeta),
\label{27}
\end{equation}
where $a, c_x, c_y$ are evaluated at $0$ and $f_i \in \MM^1$
($\MM^k$, $k \ge 0$, is the ideal of $k$-flat functions in the ring of smooth functions).

The field \eqref{24} has the invariant foliation $\zeta=\const$.
The set of integral curves of \eqref{24} passing through the origin consists of
the $\xi$-axis and the one-parameter family
\begin{equation}
\{\xi = \alf |\eta|^{\frac{3}{2}}, \ \zeta=0\}, \ \  \alf \in \bR,
\label{28}
\end{equation}
tending to the $\xi$-axis as $\alf \to \infty$.
Consider the possible mutual relationships between
the phase portrait of the germ \eqref{16} at $0$ and the $(x,y)$-plane in the space $PTM$
using the eigenvectors \eqref{26-08-2016}.
The integral curve of the field \eqref{16} corresponding to the $\xi$-axis in \eqref{24}
has a non-vertical tangential direction at $0$ (the eigenvector $e_1$),
hence its projection to the $(x,y)$-plane is a smooth geodesic.
On the contrary, the family \eqref{28} gives a family of integral curves of \eqref{16}
with vertical tangential direction at $0$ (the eigenvector $e_2$).
The projections of these curves to the $(x,y)$-plane have singularity at $0$.

To establish the character of the singularity, substitute \eqref{28} in \eqref{27}.
This yields
$$
x = 2a\alf |\eta|^{\frac{3}{2}} + \bar f_{1,\alf}(\eta), \ \ p = \eta + \bar f_{2,\alf}(\eta),
$$
where $\bar f_{i,\alf}(\eta) = f_i (\alf |\eta|^{\frac{3}{2}}, \eta,0)$.
Observe that the functions $\bar f_{1,\alf} = o\bigl(|\eta|^{\frac{3}{2}}\bigr)$ and
$\bar f_{2,\alf} = o(\eta)$ are $C^2$ and $C^1$, resp.
Denote the sign of $\eta$ by $s(\eta)$, then we have the equation
\begin{equation*}
dy = pdx = \bigl( \eta + \bar f_{2,\alf}(\eta) \bigr)
\bigl( 3a\alf |\eta|^{\frac{1}{2}} s(\eta) + \bar f'_{1,\alf}(\eta) \bigr)\, d\eta =
\bigl( 3a\alf |\eta|^{\frac{3}{2}} + g_{\alf}(\eta) \bigr)\, d\eta,
\end{equation*}
where $g_{\alf} = o\bigl(|\eta|^{\frac{3}{2}}\bigr)$ is $C^1$-smooth.
Integrating, we get
$y = \tfrac{6}{5}a\alf |\eta|^{\frac{5}{2}} s(\eta) + h_{\alf}(\eta) =
\tfrac{6}{5}a\alf \eta |\eta|^{\frac{3}{2}} + h_{\alf}(\eta)$,
where $h_{\alf} = o\bigl(|\eta|^{\frac{5}{2}}\bigr)$ is $C^2$-smooth.
The scaling $y \mapsto \tfrac{5}{3}y$, $\alf \mapsto \pm 2a\alf$ yields
\begin{equation}
x = \alf |\eta|^{\frac{3}{2}} + X_{\alf}(\eta), \ \
y = \alf \eta |\eta|^{\frac{3}{2}} + Y_{\alf}(\eta),
\label{29}
\end{equation}
where $X_{\alf} = o\bigl(|\eta|^{\frac{3}{2}}\bigr)$ and $Y_{\alf} = o\bigl(|\eta|^{\frac{5}{2}}\bigr)$
are $C^2$-smooth, and $\alf>0$ ($\alf<0$) corresponds to the domain $M_+$ ($M_-$, resp.).
The asymptotic formula \eqref{29} makes sense for all real $\alf \neq 0$.

In order to take care of the omitted  case $\alf=0$, recall that
the isotropic surface $\FF$ is an invariant surface of the field \eqref{16} (Lemma~\ref{L2})
and contains its singular points (Lemma~\ref{L4}).
Hence in the normal coordinates $(\xi,\eta,\zeta)$ the surface $\FF$ contains the $\zeta$-axis
and intersects every invariant leaf $\zeta = \const$ by a certain integral curve of \eqref{24}.
For instance, $\FF$ intersects the leaf $\zeta=0$ by an integral curve of the family \eqref{28}, which
corresponds to an isotropic geodesic passing through $0$.

On the other hand, we know that the implicit differential equation $F = ap^2 + c =0$, which
described the isotropic lines in $(M, \oF)$, has Cibrario normal form at $0$. Hence there exists
a unique isotropic geodesic passing through $0$, the semicubic parabola
\begin{equation}
x = \eta^2, \ \ y = \eta^3 N(\eta), \ \
N(0) \neq 0,
\label{30}
\end{equation}
lying in the domain $M_+$.

From the uniqueness of the isotropic geodesic passing through $0$ it follows that
the lift of \eqref{30} is the curve of the family \eqref{29} with $\alf=0$ and $Y_0(\eta)=\eta^3 N(\eta)$.
Using the representation $N(\eta) = N_1(\eta^2) + \eta N_2(\eta^2)$
and the change of variables $y \mapsto N_1(0)(y-x^2N_2(x))/N_1(x)$,
we get $N(\eta) \equiv N(0)$.
It is not hard to check that the number $\eps=N(0)$ is a mutual invariant
of the curves \eqref{29} and \eqref{30}.
\ebox

In Example \ref{Example1}, we considered $F = p^2 + c$ with $c=-x$ and $c = \alf y^2 - x$.
In both cases the isotropic direction $p_0 = 0$ is transverse to the curve $M_{0,1}$
given by the equation $x=0$ and $x = \alf y^2$, respectively, and the conditions of
Theorem~\ref{T3} hold true.

\begin{exam}
{\rm

We consider here the case $F = p^2 - x$ in more detail.
The field \eqref{16} is
\begin{equation}
-2(3x+p^2) \biggl( \frac{\pa}{\pa x} +  p \frac{\pa}{\pa y} \biggr) - 4p \frac{\pa}{\pa p}.
\label{31}
\end{equation}
It is easy to check that the isotropic surface $\FF$ given by $p^2=x$
is an invariant surface of the field \eqref{31} and the unique isotropic line passing through $0$
is given by $x = p^2$, $y = \tfrac{2}{3}p^3$.
Integrating the equation $dp/dx = 2p/(3x+p^2)$, we get the family
$x = \alf |p|^{\frac{3}{2}} + p^2$, where $\alf$ is the constant of integration,
and a single integral curve $p=0$, which gives the smooth non-isotropic geodesic $y=0$.

Integrating the relation
$
dy = pdx = p \bigl( \tfrac{3}{2}\alf |p|^{\frac{1}{2}} s(p) + 2p \bigr) dp =
\bigl( \tfrac{3}{2}\alf |p|^{\frac{3}{2}} + 2p^2 \bigr) dp,
$
we get
$y = \frac{3}{5} \alf p |p|^{\frac{3}{2}} + \frac{2}{3} p^3 + c_1$,
where $c_1$ is the second constant of integration.
The family of geodesics outgoing from $q=0$ is characterized by $c_1=0$.
The scaling $y \to \frac{5}{3}y$ brings this family to the form \eqref{25}
with $\bar X_{\alf} \equiv \bar Y_{\alf} \equiv 0$ and $\eta=p$:
\begin{equation}
x = \alf |\eta|^{\frac{3}{2}} + \eta^2, \ \ \
y = \alf \eta |\eta|^{\frac{3}{2}} + \tfrac{10}{9} \eta^3.
\label{32}
\end{equation}
For $\alf=0$, formula \eqref{32} gives the isotropic geodesic.
For $\alf \to \infty$, the curves \eqref{32} tend to the smooth geodesic $y=0$.
}
\label{Example3}
\end{exam}

\begin{remark}
\label{R25}
{\rm
Theorem~\ref{T3} shows that the extension of geodesics through the curve $M_{0,1}$ is not uniquely defined.
Indeed, all geodesics of the family \eqref{25} have the same tangential direction at $q \in M_{0,1}$ and
almost all of then have singularity of the same type at $q$. So, a curve given by formula \eqref{25} with
any $\alf \neq 0$ does not have any advantages in comparison with the curve consisting of two bows
\eqref{25} with $\alf_1$ if $(x,y) \in M_+$ and $\alf_2$ if $(x,y) \in M_-$.
}
\end{remark}

\section{Surfaces in Berwald--Moor spaces}\label{sect3}

Consider the space $\bR^n$, $n\ge 3$, with the coordinates $(x_1, \ldots, x_n)$ equipped with
{\it Berwald--Moor metric}
$ds = (dx_1 \, \cdots \, dx_n)^{\frac{1}{n}}$,
and a smooth two-dimensional surface $M \subset \bR^n$ parametrized by
$x_i = f_i(x,y)$, $i=1, \ldots, n$.
The Berwald--Moor metric of the ambient space defines two-dimensional pseudo-Finsler space $(M, \oF)$
with the metric function $\of = \oF^{\frac{1}{n}}$, where
\begin{equation}
\oF(x,y;\dx,\dy)= \prod_{i=1}^n \bigl(f_{ix}(x,y)\dx + f_{iy}(x,y)\dy \bigr), \ \
f_{ix} = \frac{\pa f_i}{\pa x}, \ \ f_{iy} = \frac{\pa f_i}{\pa y},
\label{33}
\end{equation}
and $n$ families of isotropic lines
\begin{equation}
f_i(x,y) = \const, \ \ i = 1, \ldots, n.
\label{34}
\end{equation}

Given $q \in M$, the isotropic direction $p$ is called {\it simple} ({\it double} or {\it multiple})
if there exist {\it only one} ({\it only two} or {\it more than one}, resp.) isotropic lines \eqref{34}
passing through $q$ with given direction $p$.
By Lemma~\ref{L4}, singularities of the geodesic flow occur at the points $q \in M$ that
have at least one multiple isotropic direction.

\begin{remark}
\label{R3}
{\rm
In the case $n=3$ we have a cubic pseudo-Finsler space $(M,\oF)$.
But unlike Section~\ref{sect2.2}, the function $\oF$ given by \eqref{33} is not generic.
The corresponding cubic polynomial $F(x,y;p)$ at every point $q \in M$ has $n$ real roots
(taking into account the multiplicity and including the root $p=\infty$),
and $M_- = \emptyset$. Hence singularities of geodesic appear only at the points
where at least two of three isotropic lines \eqref{34} are tangent.
Here the stratum $M_{0,1}$ consists of the points where two isotropic lines are tangent
(the {\it double} isotropic direction) and
the third one is transversal to them (the {\it simple} isotropic direction).
}
\end{remark}

From now on, we assume that the functions $f_i(x,y)$ have non-degenerate differentials and
every point $q \in M$ may have simple or double isotropic directions only
(the number of double isotropic directions can vary from $0$ to $[\frac{n}{2}]$).
Moreover, assume that the tangency of isotropic lines with double isotropic directions has first order.
Consider geodesics passing through a point $q$ with a double isotropic direction $p_0$
satisfying the above conditions.

Without loss of generality assume that $q=0$ (the origin in the $(x,y)$-plane) and $p_0$
corresponds to the isotropic lines $f_1(x,y) = 0$ and $f_2(x,y) = 0$, where $f_{2y}(0) \neq 0$.
Making the change of variable $y \mapsto f_2(x,y)$, we transform the metric function \eqref{33} into
a similar one with $f_2(x,y)=y$ and $f_{1x}(0)=0$, $f_{1y}(0) \neq 0$, $f_{1xx}(0) \neq 0$.
The double isotropic direction $p_0$ becomes $p=0$ and, moreover,
in a neighborhood of $0$, $p=0$ is the double isotropic direction at all points
$q \in M_{0,1} = \{ f_{1x}(x,y)=0 \}$.

By the condition $f_{1xx}(0) \neq 0$, $M_{0,1}$ is a smooth curve transversal to the $x$-axis.
Making the change of variable $x \mapsto f_{1x}(x,y)$, we transform $M_{0,1}$ into $x=0$ and
the metric function \eqref{33} into a similar one with
$f_{1x}=xa(x,y)$, $f_{1y}=b(x,y)$, $f_2(x,y)=y$,
where $a,b$ are smooth functions non-vanishing at $0$.
So, we get
\begin{equation}
F(x,y;p) = p(ax+bp)G, \ \ \ G(x,y;p) = \prod_{i=3}^n (f_{ix}+f_{iy}p), \ \ \, G(0,0;0) \neq 0.
\label{35}
\end{equation}

Substituting \eqref{35} in \eqref{13}, we get
\begin{equation}
\label{36}
\begin{aligned}
\Del &= \bigl( (1-n)(ax)^2 + 2(2-n)abxp + 2(2-n)(bp)^2 \bigr)\, G^2 + \Del_0,   \\
   P &= ap \bigl( (n-2)bp - ax + P_0 \bigr) \, G^2,   \\
\end{aligned}
\end{equation}
where $\Del_0 \in \langle x^3, x^2p, xp^2, p^3 \rangle$ and $P_0 \in \langle x^2, xp, p^2 \rangle$
(both ideals are in the ring of smooth functions on $x,y,p$).
Formula \eqref{36} shows that all components of the field \eqref{16} vanish on the line $\{x=p=0\}$
and the spectrum of \eqref{16} at any point of this line has three zero eigenvalues.
This does not allow to establish a normal form similarly to Theorem~\ref{T3}.

To overcome this problem, consider the blowing up
\begin{equation}
\FFF \colon (x,y,u) \mapsto (x,y,p), \ \ p = xu, \ \ u \in \bR P^1 = \bR \cup \infty.
\label{37}
\end{equation}
The mapping $\FFF$ is one-to-one except on the plane $\Pi = \{(x,y,u) \colon x=0\}$,
whose image is the line $\FFF (\Pi) = \{ (x,y,p) \colon x = p = 0 \}$.
The mapping $\FFF$ is a local diffeomorphism at all points of the $(x,y,u)$-space except $\Pi$.
It has an inverse defined on $\mathbb R^3\setminus \FFF (\Pi)$ given  by
$$
\FFF^{-1}(x,y,p) = \Bigl(x, y, \frac{p}{x} \Bigr).
$$

Observe that there is no geodesic that coincides with the line $\FFF (\Pi)$.
A straightforward calculation shows that the field \eqref{16} corresponds to
a smooth field in the $(x,y,u)$-space (away of $\Pi$) that,
after dividing by the common factor $xG^2$, is
\begin{equation}
x \bigl( \AA + \ldots \bigr) \biggl( \frac{\pa}{\pa x} +  xu \frac{\pa}{\pa y} \biggr) +
u(n-2) \bigl( 2(bu)^2 + 3abu + a^2 + \ldots \bigr) \frac{\pa}{\pa u},
\label{38}
\end{equation}
where
$$
\AA(x,y,u) = (1-n)a^2 + 2(2-n) (abu + (bu)^2),
$$
here and below the dots mean terms that belong to the ideal $\langle x \rangle$.

\begin{remark}
\label{R4}
{\rm
There exists $\eps>0$ such that $\AA(x,y,u) \neq 0$ for all $u$ if $|x|+|y| < \eps$.
Indeed, consider $\AA(x,y,u)$ as a quadratic polynomial on $u$ with the discriminant
$$
D = (2-n)^2(ab)^2 - 2(1-n)(2-n)(ab)^2 = n(2-n)(ab)^2,
$$
which is strictly negative if $x,y$ are sufficiently close to zero.
}
\end{remark}

Dividing the field \eqref{38} by $(\AA + \ldots)$, we get
\begin{equation}
x \biggl( \frac{\pa}{\pa x} +  xu \frac{\pa}{\pa y} \biggr) + u(n-2) (U+\ldots) \frac{\pa}{\pa u},  \ \ \
U(x,y,u) = \frac{2(bu)^2 + 3abu + a^2}{\AA(x,y,u)}.
\label{39}
\end{equation}

\begin{remark}
\label{R5}
{\rm
The plane $x=0$ is invariant for the fields \eqref{38} and \eqref{39}.
Moreover, it is filled with {\it vertical} (i.e., parallel to the $u$-direction) straight integral lines
of these fields, whose projections to the $(x,y)$-plane along the $u$-axis are points on the $y$-axis.
}
\end{remark}

\begin{lemma}
\label{L6}
Geodesics can pass through a point $q \in M$ lying on the $y$-axis with the direction $p=0$
only with the following {\it admissible values} $u:$
\begin{equation}
u_0=0, \ \ u_1 = -a/2b, \ \ u_2 = -a/b.
\label{40}
\end{equation}
\end{lemma}

{\bf Proof}.
By the standard existence and uniqueness theorem, for every point $(x,y,u)$
such that $x=0$ and $U(x,y,u) \neq 0$ there exists a unique integral curve of the field \eqref{39}
passing through this point. By Remark~\ref{R5}, it is a vertical straight line,
whose projection to the $(x,y)$-plane is a point on the $y$-axis.
Hence geodesics can pass through a point $q \in M$ lying on the $y$-axis with the direction $p=0$
only with $u = 0$ or $u$ such that $2(bu)^2 + 3abu + a^2 = 0$.
This gives the three values in \eqref{40}.
\ebox

\begin{lemma}
\label{L7}
The set of singular points of the field \eqref{39} consists of three
mutually disjoint curves
$$
W^c_i = \{ (x,y,u) \colon x=0, \, u=u_i(y)\}, \ \ \, i=0,1,2.
$$
On every curve $W^c_i$, the linear part of the field \eqref{39} has the constant spectrum $(1,\la, 0)$, where
$\la = \tfrac{n-2}{n}$ if $i=1$ or $\la = \tfrac{n-2}{1-n}$ if $i \in \{0,2\}$.
In both cases, $\pa_u$ is the eigenvector with~$\la$.
\end{lemma}

{\bf Proof}.
The first statement is trivial.
All other statements are by direct calculations.
\ebox

\begin{theorem}
\label{T4}
Suppose that the functions $f_i(x,y)$ have non-degenerate differentials
and $p_0$ is a double isotropic direction at $q \in M$ such that
the corresponding isotropic lines have first order of tangency at $q$.
Then the field \eqref{39} at its singular points $(q; u_i)$ has local orbital normal forms
indicated in Table~\ref{Tab1}
and to $p_0$ corresponds a one-parameter family of geodesics outgoing from $q$.
There exist smooth local coordinates centered at $q$
such that this family consists of $C^2$-smooth non-isotropic geodesics
\begin{equation}
y = x^2 + Y(x,\alf|x|^{\la}), \ \ \, Y(x,\alf|x|^{\la})=o(x^2), \ \ \,
\la = \tfrac{n-2}{n}, \ \ \, \alf \in \bR,
\label{41}
\end{equation}
where $Y(\cdot,\cdot)$ is a $C^{\infty}$-smooth function,
together with two $C^{\infty}$-smooth isotropic geodesics
\begin{equation}
y = 0 \, \ \ \textrm{and} \ \ \,  y = 2x^2 + o(x^2).
\label{42}
\end{equation}
In a neighborhood of $q$, every geodesic of the family \eqref{41}
belongs to the curvilinear tongue-like sector bounded by the curves \eqref{42}
as it is presented in Figure~\ref{fig4} (right).
\end{theorem}

\begin{table}[htb]
\begin{center}
\begin{tabular}{|c|c|c|}
\hline
$\phantom{\Bigl|}$\! & \multicolumn{2}{c|}{\textrm{Orbital normal form}} \\
\hline
{} $\phantom{\Bigl|}$\! & topological & $C^{\infty}$-smooth  \\
\hline
$\phantom{\Bigl|}\! W^c_1$ &
$\xi \frac{\pa}{\pa \xi} + \eta \frac{\pa}{\pa \eta} + 0\frac{\pa}{\pa \zeta}$ &
$\xi \frac{\pa}{\pa \xi} + \la \eta \frac{\pa}{\pa \eta} + 0\frac{\pa}{\pa \zeta}$, where $\la = \tfrac{n-2}{n}$ \\
\hline
$\phantom{\Bigl|}\! W^c_{0}$ & {} &
$\xi (n-1+\Phi_1(\rho,\zeta))\frac{\pa}{\pa \xi} - \eta (n-2+\Phi_2(\rho,\zeta)) \frac{\pa}{\pa \eta} +
\rho \Psi(\rho,\zeta) \frac{\pa}{\pa \zeta}$,
\\
$\phantom{\Bigl|}\!$ and $\phantom{\Bigl|}$\! &
$\xi \frac{\pa}{\pa \xi} - \eta \frac{\pa}{\pa \eta} + 0\frac{\pa}{\pa \zeta}$ &
\textrm{where} \ $\rho = \xi^{n-2}\eta^{n-1}$; \\
$\phantom{\Bigl|}\! W^c_{2}$ &
{} &
$(n-1)\xi \frac{\pa}{\pa \xi} - (n-2) \eta \frac{\pa}{\pa \eta} + \rho \frac{\pa}{\pa \zeta}$, \
\textrm{if} \ $\Psi(0,0) \neq 0$. \\
\hline
\end{tabular}
\caption{Local orbital normal forms of the field \eqref{39}.}
\label{Tab1}
\end{center}
\end{table}

{\bf Proof}.
Choose local coordinates so that $q \in M$ is the origin and consider the field \eqref{39} in a neighborhood
of its singular points $(0,0,u_i)$, $i=0,1,2$, where $u_i$ are given by formula \eqref{40}.
By Lemma~\ref{L7}, in all singular points the condition $\Real \la_{1,2} \neq 0$ holds,
and every curve $W^c_i$, $i=0,1,2$, is the center manifolds of this field.
Moreover, there exist also 2-dimensional unstable manifold if $i=1$ and the pair of 1-dimensional
stable and unstable manifolds if $i=0,2$.
Hence all topological normal forms in Table~\ref{Tab1} trivially follow
from the reduction principle \cite{A-Geom, AI, HPS}.

Indeed, the reduction principle asserts that the germ \eqref{39} is orbitally topologically equivalent
to the direct product of the standard 2-dimensional node (if $i=1$) or saddle (if $i=0,2$)
and the restriction of the field to the center manifold $W^c_i$.
Since the restriction of the field \eqref{39} to every center manifold $W^c_i$, $i=0,1,2$,
is identically zero, this gives us the topological normal forms in Table~\ref{Tab1}.

Establish now the smooth normal forms in the cases $i=1$ and $i=0,2$ separately.

{\it The case $i=1$}.
By Lemma~\ref{L7}, the linear part of the field \eqref{39} at any point on $W^c_1$ has
spectrum $\bigl(1, \la, 0\bigr)$ with $\la = \tfrac{n-2}{n}$.
Then Theorem~5.5 in \cite{RemizovTari} asserts that the germ \eqref{39} at any point on $W^c_1$
is orbitally $C^{\infty}$-smoothly equivalent to
\begin{equation}
(\xi + \phi(\zeta)\eta^{1/\la}) \frac{\pa}{\pa \xi} + \la \eta \frac{\pa}{\pa \eta} + 0 \frac{\pa}{\pa \zeta},
\label{43}
\end{equation}
where $\phi(\zeta) \equiv 0$ if the number $1/\la$ is not integer (non-resonant case).

Assume $1/\la$ is integer and prove that $\phi(\zeta) \equiv 0$ iff
for every point $\om_* \in W^c_1$ the field \eqref{39} has a $C^{\infty}$-smooth integral curve
passing through $\om_*$ with the vertical tangential direction $\pa_u$. By Remark~\ref{R5},
such integral curve exists (the vertical straight lines),
hence this establishes the equality $\phi(\zeta) \equiv 0$
in the remaining cases $1/\la=3$ ($n=3$) and $1/\la=2$ ($n=4$).

For this, note that the field \eqref{43} has the invariant foliation $\zeta = \const$.
Every invariant leaf contains a single integral curve $\eta=0$
corresponding to eigendirection with the eigenvalue $1$
and one-parameter family of integral curves
\begin{equation}
\xi = \eta^{1/\la} (\alf + \phi(\zeta)\ln |\eta|), \ \ \alf \in \bR,
\label{44}
\end{equation}
corresponding to eigendirection with the eigenvalue $\la$.
All curves \eqref{44} are $C^{1/\la-1}$-smooth (but not $C^{1/\la}$-smooth at zero) if $\phi(\zeta) \neq 0$
and $C^{\infty}$-smooth if $\phi(\zeta)=0$.
Without loss of generality, assume that the point $\om_* \in W^c_1$ in the $(x,y,u)$-space corresponds
to $(0,0,\zeta_*)$ in the $(\xi,\eta,\zeta)$-space.
The equality $\phi(\zeta_*)=0$ is equivalent to the existence
of at least one $C^{\infty}$-smooth integral curve of the field \eqref{43} with tangential direction $\pa_{\eta}$
lying on the invariant leaf $\{\zeta = \zeta_*\}$.
To complete the proof, remark that the eigendirection $\pa_{\eta}$ of \eqref{43}
corresponds to the eigendirection $\pa_{u}$ of \eqref{39}.

{\it The cases $i=0,2$}.
By Lemma~\ref{L7}, the linear part of the field \eqref{39} at all points on the curves $W^c_0$ and $W^c_2$
has spectrum $(\la_1,\la_2,0)$, where $\la_1=1$ and $\la_2=\tfrac{n-2}{1-n}$.
This gives the resonance
\begin{equation}
\mu \la_1 + \nu \la_2=0
\label{res}
\end{equation}
with the resonant monomial $\rho = \xi^{\mu}\eta^{\nu}$,
where we set $\mu = n-2$ and $\nu = n-1$.
Everything that we say below is true as well for arbitrary relatively prime $\mu, \nu \in \bN$.

The resonance \eqref{res} does not allow to get a normal form with one identically zero component
(as we have in the case $i=1$) even in the finite-smooth category, see the discussion in \cite{GR} (Section~3.2).
Moreover, \eqref{res} generates two infinite series of resonances
\begin{equation*}
(1+l\mu)\la_1+l\nu\la_2=\la_1,  \quad  l\mu\la_1+(1+l\nu)\la_2=\la_2,  \quad l=1,2, \ldots,
\end{equation*}
and consequently, infinite number of resonant monomials in the corresponding (orbital)
Poincar\'e--Dulac normal form:
\begin{equation}
\xi (\nu+\Phi_1(\rho,\zeta))\frac{\pa}{\pa \xi} - \eta (\mu+\Phi_2(\rho,\zeta)) \frac{\pa}{\pa \eta} +
\rho \Psi(\rho,\zeta) \frac{\pa}{\pa \zeta},
\label{27-08-2016}
\end{equation}
see e.g. \cite{GR} (Section~3.2) or \cite{Rem06} (Section~5).

Moreover, if in addition, $\Psi \neq 0$ at a point $\om_*$, the germ \eqref{27-08-2016} at $\om_*$
is smoothly orbitally equivalent to
\begin{equation}
\nu\xi \frac{\pa}{\pa \xi} - \mu \eta \frac{\pa}{\pa \eta} + \rho \frac{\pa}{\pa \zeta}.
\label{27-08-2016-0}
\end{equation}
The normal form \eqref{27-08-2016-0} was firstly established by Roussarie  \cite{Rouss}
in the partial case $\mu=\nu=1$ in $C^{\infty}$-smooth category.
For arbitrary integers $\mu, \nu$, the proof (in finite-smooth category) can be found in \cite{Rem06} (Section~5).
Combining the methods from \cite{Rouss} and \cite{Rem06}, one can
establish the normal form \eqref{27-08-2016-0} with arbitrary $\mu, \nu$
in $C^{\infty}$-smooth category also.

{\it Completion of the proof}.
Integral curves of the field \eqref{39} passing through $(0,0,u_1)$ correspond to integral curves
of the field $\xi \frac{\pa}{\pa \xi} + \la \eta \frac{\pa}{\pa \eta}$ lying on the invariant leaf $\{\zeta = 0\}$:
a single curve that coincides with the $\eta$-axis
and one-parameter family $\{\eta= \alf|\xi|^{\la}, \ \zeta=0\}$, $\alf \in \bR$.
Comparing the germ \eqref{39} at $(0,0,u_1)$ with its normal form
$\xi \frac{\pa}{\pa \xi} + \la \eta \frac{\pa}{\pa \eta}$,
one can see that the conjugating diffeomorphism $(x,y,u) \mapsto (\xi, \eta, \zeta)$
can be chosen in the form
\begin{equation}
x = \xi, \ \
u = u_1 + \eta + c_1 \xi + c_2 \zeta + \phi(\xi,\eta,\zeta), \ \
y = \zeta + \psi(\xi,\eta,\zeta),  \ \ \,
\phi,\psi \in \MM^2,
\label{45}
\end{equation}
where $\MM^k$, $k \ge 0$, is the ideal of $k$-flat functions in the ring of smooth functions.
Substituting $\xi=\zeta=0$ in \eqref{45} and taking into account $p=ux$, we get $x=p=0$.
Hence the $\eta$-axis does not correspond to a geodesic.

Substituting $\eta= \alf|\xi|^{\la}$ and $\zeta=0$ in \eqref{45} and taking into account $x=\xi$ and $p=ux$,
we get $p = x (u_1 + f(x,\alf|x|^{\la}))$ with a certain smooth function $f \in \MM^0$.
This gives the relation $dy = p dx = x (u_1 + f(x,\alf|x|^{\la})) dx$,
where  $xf(x,\alf|x|^{\la}) = o(x)$ is a $C^1$-smooth function.
Integrating, we get $y = \frac{u_1}{2} x^2 + Y(x,\alf|x|^{\la})$.
Here $Y(\cdot,\cdot)$ is a smooth function and $Y(x,\alf|x|^{\la})=o(x^2)$ is $C^2$-smooth.
After the scaling $y \to 2y/{u_1}$, we get the family \eqref{41}.

The topological and smooth orbital normal forms in Table \ref{Tab1} show that the field \eqref{39}
has only two integral curves passing through its singular point $(0,0,u_i)$, where $i=0$ or $2$.
Moreover, one of these integral curves is straight vertical line, whose projection to the $(x,y)$-plane
is a point (see Remark~\ref{R5} and Lemma~\ref{L7}). Another integral curve has non-vertical tangential
direction at $(0,0,u_i)$, hence its projection to the $(x,y)$-plane is regular.

Thus every of the admissible values $u_0$ and $u_2$ gives a smooth geodesic passing through the point $q$
with tangential direction $p=0$. It is not hard to see that these geodesics are isotropic lines, which
are solutions of differential equations $p=0$ and $ax+bp=0$, respectively (see formula~\eqref{35}).
Taking into account \eqref{40}, after the scaling $y \to 2y/{u_1}$ we get \eqref{42}.
\hfill\ebox

\begin{remark}
\label{R6}
{\rm
The normal forms
$3\xi \frac{\pa}{\pa \xi} + 2\eta \frac{\pa}{\pa \eta}$ and
$\xi \frac{\pa}{\pa \xi} + \frac{n-2}{n} \eta \frac{\pa}{\pa \eta}$
in Theorems \ref{T3}, \ref{T4}
are valid also in the analytic category; see e.g. \cite{Voronin}.
Therefore, in the analytic case, formulae \eqref{25} and \eqref{41} present Puiseux series for geodesics.
}
\end{remark}

\begin{exam}
{\rm

Consider geodesics on the surface $z=y-2x^2$ in the Berwald--Moor space
$(x,y,z)$ with the metric $ds = (dx\,dy\,dz)^{\frac{1}{3}}$.
This yields
\begin{equation}
F(x,y;p) = p(p-4x), \ \  \Del = -2(p^2-4xp +16x^2), \ \ P = -4p(p+4x),
\label{Vavilon}
\end{equation}
and equation of geodesics \eqref{12} reads
\begin{equation}
\frac{dp}{dx} = \frac{2p(p+4x)}{p^2-4xp +16x^2}.
\label{46}
\end{equation}

The isotropic lines are solution of the differential equation $p(p-4x)=0$.
It gives two families of isotropic lines $y = \const$ and $y = 2x^2 + \const$,
which have the first order tangency on the line $x=0$.
Substituting them into \eqref{46}, one can see that they are geodesics.

Consider the geodesics outgoing from the point $q=0$ with the double isotropic directions $p_0=0$.
(Recall that for every $p \neq 0$ there exists a unique geodesic passing through $q$ with tangential direction $p$,
we exclude such geodesics from further consideration.)
The isotropic geodesics $y = 0$ and $y = 2x^2$ (formula \eqref{42}) separate the $(x,y)$-plane
into four parts: the upper domain $y > 2x^2$, the semiplane $y<0$ and two tongue-like sectors between them.
See Fig.~\ref{fig4} (right), the isotropic geodesics $y = 0$ and $y = 2x^2$ are depicted with dashed lines.

Theorem \ref{T4} claims that there exists a one-parameter family of geodesics outgoing from $q$
with the double isotropic directions $p_0=0$
into the tongue-like sectors (non-isotropic family \eqref{41}) and there are no geodesics outgoing from $q$
with the double isotropic directions $p_0=0$
into two remaining parts of the plane.
Geodesics of the family \eqref{41} correspond the admissible value $u_1=2$
(compare formulae \eqref{35}, \eqref{40} and \eqref{Vavilon})
and they can be presented as the Puiseux series
$$
y = t^6 + 3t^6 \sum_{i \ge 4} \frac{a_i}{i+3} t^{i-3}, \ \ \,
p = \frac{dy}{dx} = \frac{1}{3t^2} \frac{dy}{dt} = 2t^3 + \sum_{i \ge 4} a_i t^i,  \ \ \,
\textrm{where} \ \ \,  x=t^3.
$$
Substituting the above expression for $p$ in \eqref{46}, we obtain
recurrence relations for the unknown coefficients $a_i$.

Namely, $a_i=0$ for all odd $i$ (this also follows from the fact that the surface $z=y-2x^2$
is symmetric with respect to the plane $x=0$). For even $i$ we have
$24a_6 + 4a_4^3 = 0$, $48a_8 + 14a_4^2 a_6 = 0$, $72a_{10} + 16(a_4 a_6^2 + a_4^2 a_8) = 0$, etc.
In general,
\begin{equation}
12(2i-4) a_{2i} + b_{2i} = 0, \ \ \, i=3,4,5, \ldots,
\label{47}
\end{equation}
where $b_{2i}$ is a polynomial on $a_{2j}$ with $j<i$ with zero free term.
This shows that the coefficient $a_4$ is arbitrary,
and all $a_{2i}$ with $i \ge 3$ are uniquely defined by equations \eqref{47}.
This gives the one-parameter non-isotropic family \eqref{41}.
In particular, $a_4=0$ gives $a_{2i}=0$ for all $i \ge 3$, and the corresponding solution $y=x^2$
presents the unique $C^{\infty}$-smooth geodesic of non-isotropic family \eqref{41}
(the corresponding value of the parameter is $\alf=0$).
}
\label{Example4}
\end{exam}





\normalsize

\begin{thebibliography}{99}


\bibitem{AM}
\emph{Antonelli P. L., Miron R.},
Lagrange and Finsler Geometry. Applications to Physics and Biology.
Kluwer Academic Publishers, 1996.

\bibitem{A-Geom}
\emph{Arnol'd V. I.},
Geometrical methods in the theory of ordinary differential equations.
Springer-Verlag, New York, 1988.

\bibitem{AGiv}
\emph{Arnol'd V. I., Givental' A. B.},
Symplectic geometry.
Encyclopedia of Mathematical Sciences, Dynamical Systems IV (Springer, Berlin, 1985), pp.~5--131.

\bibitem{AI}
\emph{Arnol'd V. I., Il'yashenko Yu. S.},
{Ordinary differential equations, Dynamical systems~I}. Encyclopaedia Math. Sci., vol. 1,
Springer-Verlag 1988.

\bibitem{Asanov}
\emph{Asanov G. S.},
Finsler Geometry, Reltivity and Gauge Theories.
Reidel, Dordrecht, 1985.

\bibitem{Balan}
\emph{Balan V., Neagu M.},
Jet single-time Lagrange geometry and its applications. John Wiley \& Sons, Inc., Hoboken, NJ, 2011.

\bibitem{BCS}
\emph{Bao D., Chern S.-S., Shen Z.},
An Introduction to Riemann-Finsler Geometry.
Graduate Texts in Mathematics, 200. Springer-Verlag, New York, 2000.

\bibitem{GR}
\emph{Ghezzi R., Remizov A. O.},
On a class of vector fields with discontinuities of divide-by-zero type and its applications to
geodesics in singular metrics. Journal of Dynamical and Control Systems, 18:1 (2012), pp.~135--158.

\bibitem{HPS}
\emph{Hirsch M. W., Pugh C. C., Shub M.},
{Invariant manifolds.}
Lecture Notes in Mathematics, Vol. 583. Springer-Verlag, Berlin-New York, 1977.

\bibitem{KhT}
\emph{Khesin B., Tabachnikov S.},
Pseudo-Riemannian geodesics and billiards.
Advances in Math. 221 (2009), pp.~1364--1396.

\bibitem{KosKri}
\emph{Kossowski M., Kriele M.},
Smooth and discontinuous signature type change in general relativity.
Class. Quantum Grav. 10 (1993), pp.~2363--2371.

\bibitem{Martinet}
\emph{Martinet J.},
Sur les singularit\'es des formes diff\'erentielles,
Ann. Inst. Fourier, 1970, 20, no\,1, pp.~95--178.

\bibitem{Matsumoto95}
\emph{Matsumoto M.},
Two-dimensional Finsler spaces whose geodesics constitute a family of special conic sections.
J. Math. Kyoto Univ. 35 (1995), no.~3, pp.~357--376.

\bibitem{Matsumoto-Shimada}
\emph{Matsumoto M., Shimada H.},
On Finsler spaces with 1-form metric. II. Berwald-Moor's metric $L = (y^1y^2 \cdots y^n)^{1/n}$.
Tensor (N.S.) 32 (1978), no.~3, pp.~275--278.

\bibitem{MHV}
\emph{Mike\v{s} J., Hinterleitner I., Van\v{z}urov\'a A.},
One remark on variational properties of geodesics in pseudoriemannian and generalized Finsler spaces.
In: Geometry, integrability and quantization, Softex, Sofia, 2008, pp.~261--264.

\bibitem{PR}
\emph{Pavlova N. G., Remizov A. O.},
Geodesics on hypersurfaces in the Minkowski space: singularities of signature change.
Russian Math. Surveys, 2011, 66:6 (402), pp.~193--194.

\bibitem{Rem06}
\emph{Remizov A. O.},
Multidimensional Poincar\'e construction and singularities of lifted fields for implicit differential equations.
J. Math. Sci. (N.Y.) 151:6 (2008), pp.~3561--3602.

\bibitem{Rem08}
\emph{Remizov A. O.},
Codimension-two singularities in 3D affine control systems with a scalar control.
{Mat. Sb.}, 199:4 (2008), 143--158.

\bibitem{Rem-Pseudo}
\emph{Remizov A. O.},
Geodesics on 2-surfaces with pseudo-Riemannian metric: singularities of changes of signature.
{Mat. Sb.}, 200:3 (2009), 75--94.

\bibitem{Rem15}
\emph{Remizov A. O.},
On the local and global properties of geodesics in pseudo-Riemannian metrics.
{Differential Geometry and its Applications}, 39 (2015), 36--58.

\bibitem{RemizovTari}
\emph{Remizov A. O., Tari F.},
Singularities of the geodesic flow on surfaces with pseudo-Riemannian metrics.
{Geometriae Dedicata}, 185 (2016), no.~1, 131--153.

\bibitem{Rouss}
\emph{Roussarie R.},
Mod\`eles locaux de champs et de formes.
Asterisque, vol.\,30 (1975), pp.~1--181.

\bibitem{Rund}
\emph{Rund H.},
The differential geometry of Finsler spaces.
Die Grundlehren der Mathematischen Wissenschaften, Bd. 101.
Springer-Verlag, Berlin-G\"ottingen-Heidelberg, 1959, xiii+284 pp.

\bibitem{ZhitSo}
\emph{Sotomayor J., Zhitomirskii M.},
Impasse singularities of differential systems of the form $A(x)x'=F(x)$.
J. Diff. Equations 169 (2001), no.\,2, pp.~567--587.

\bibitem{Voronin}
\emph{Voronin S. M.},
The Darboux--Whitney theorem and related questions.
In: Nonlinear Stokes phenomenon (Yu.S.~Ilyashenko, ed.).
Adv. Sov. Math.~14, Providence (1993), pp.~139--233.

\end{thebibliography}
\end{document}